[1994/12/01]
\documentclass[manuscript]{aomart}
\usepackage{graphicx}
\usepackage{amsmath,amsfonts,amssymb,amscd,amsthm,xspace}
\usepackage[utf8]{inputenc}
\chardef\bslash=`\\ 

\newtheorem[{}\it]{thm}{Theorem}[section]

\newtheorem{lem}[thm]{Lemma}
\newtheorem{prop}[thm]{Proposition}

\theoremstyle{definition}

\newtheorem*[{}\it]{notation}{Notation}

\newcommand{\thmref}[1]{Theorem~\ref{#1}}

\newcommand{\lemref}[1]{Lemma~\ref{#1}}
\newcommand{\pref}[1]{Proposition~\ref{#1}}
\newcommand{\eref}[1]{Equation~\ref{#1}}

\newcommand{\eval}[2][\right]{\relax
  \ifx#1\right\relax \left.\fi#2#1\rvert}

\title[Brownian motion and Symmetrization]{Brownian motion and Symmetrization}
\author{Tomas Kojar}
\address{University of Toronto, Canada}
\fulladdress{Graduate Department of Mathematics,University of Toronto, M5S 2E4 , ON , Canada}
\email{Tomas.Kojar@mail.utoronto.ca}
\givenname{Tomas}
\surname{Kojar}


\publicationyear{2015}
\papernumber{12}
\startpage{1}
\endpage{}
\begin{document}
\maketitle
\begin{abstract}
  In this survey we explore the salient connections made between Brownian motion, symmetrization and complex analysis in the last 60 years starting with Kakutani's paper (1944) equating harmonic measure and exit probability. To exemplify these connections we will survey the techniques used in the literature to prove isoperimetric results for exit probabilities and Riesz capacities.

\end{abstract}

\tableofcontents

\section{Introduction}

The classical isoperimetric problem asks: Given all shapes of a given area, which of them has the minimal perimeter. The conjectured answer was the disk and Steiner in 1838 showed this to be true using the Steiner symmetrization method (described below). From this many other isoperimetric problems sprung. We will study the following three areas: Rayleigh's conjecture (1877) that the first eigenvalue of the Dirichlet problem is minimized for the ball was proved independently by G. Faber and E. Krahn. P$\acute{o}$lya and G. Szeg$\ddot{o}$ (1951) proved that for fixed volume, the ball has the minimum electrostatic capacity.  Finally, many isoperimetric results for the green function and harmonic measure proved by Luttinger, Baernstein and many other authors mentioned throughout this survey. \\
The first person to describe the mathematics behind Brownian motion was the Danish astronomer Thorvald Thiele in 1880, and later, in 1900, Louis Bachelier a French mathematician, wrote his PhD thesis on the "Theory of Speculation", which was the first ever mathematical analysis of the stock and option markets. Bachelier's work also provided a mathematical account of Brownian Motion. Einstein and Smoluchowski (1906) realised that movements of Brownian particles were caused by collisions with molecules of the solvent. These molecules move erratically in display of their thermal energy, of which the temperature is a certain measure.  Today this explanation may seem to be trivial, but a hundred years ago the atomistic hypothesis was not commonly accepted. Finally, Wiener took a great interest in the mathematical theory of Brownian motion proving many results now widely known such as the non-differentiability of the paths. Consequently the one-dimensional Brownian motion was named the Wiener process. It is the best known of the L$\acute{e}$vy processes, c$\grave{a}$dl$\grave{a}$g stochastic processes with stationary statistically independent increments, and occurs frequently in pure and applied mathematics, physics and economics. \\
The connection of these two areas crosses through the complex space and it started with Kakutani's result (1944) (proved below) of representing solutions to the Dirichlet problem as expected value of Brownian motion at the boundary. From then on, many complex analysis and pdes' objects were phrased probabilistically. Green function $G(x,y)$ as the expected total number of visits to y, starting from x. Electrostatic capacity of solid A as the total heat A can absorb proved by F. Spitzer (1964)\cite{Spitzer_1964}. The first eigenvalue of set A as the asymptotic probability of a Brownian particle entering A proved by Kac (1951)\cite{Kac_1966}. These and other connections will be surveyed in this thesis. In this section, we describe the symmetrization methods we will need and prove Kakutani's result on the Dirichlet problem.

\subsection{Symmetrization methods}
Let $\Omega\subset \mathbb{R}^{n}$ be measurable, then we denote by $\Omega^{*}$ the \textit{symmetrized } version of $\Omega$ i.e. a ball $\Omega^{*}:=B_{r}(0)\subset\mathbb{R}^{n}$ such that $vol(\Omega^{*})=vol(\Omega)$. We denote by $f^{*}$ the \textit{symmetric decreasing rearrangement} of nonnegative measurable function f and define it as $f^{*}(x):=\int_{0}^{\infty}1_{\{f(x)>t\}^{*}}(t)dt$. The following methods have been proved to transform $\Omega$ to $\Omega^{*}$ i.e. given a sequence of symmetrization transformations $\{T_{k}\}$ we have $\lim\limits_{k\to \infty}d_{Ha}(\Omega^{*}, T_{k}(K) )=0$, where $d_{Ha}$ is the Hausdorff distance \cite{Burchard_2009}. 

\begin{figure}[h]
\centering
\includegraphics[scale=0.3]{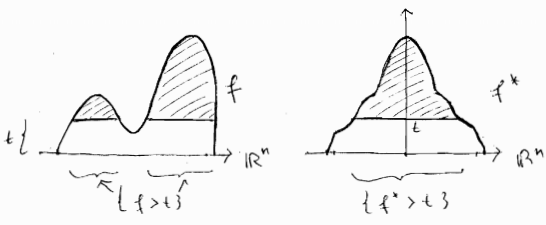}
\caption{Symmetric decreasing rearrangement of f}
\end{figure}

\subsubsection{Steiner symmetrization}
Steiner symmetrization was introduced by Steiner (1838) to solve the isoperimetric theorem stated above. Let $H^{n-1}\subset\mathbb{R}^{n}$ be a hyperplane through the origin. Rotate space so that $H^{n-1}$ is the $x_{n}=0$ hyperplane. For each $x\in H$ let the perpendicular line through $x\in H $ be $L_{x}=\{x+ye_{n}:y\in \mathbb{R}\}$. Then by replacing each $\Omega\cap L_{x}$ by a line centered at H and with length $|\Omega\cap L_{x}|$ we obtain the \textit{Steiner symmetrized} version. We denote by $St(f)$ the Steiner symmetrization wrt to $x_{n}=0$ hyperplane of nonnegative measurable function $f:\mathbb{R}^{d}\to \mathbb{R}$ and define it as $St(f)(x_{1},...,x_{n}):=f^{*}(x_{1},...,x_{n})$ for fixed $x_{1},...,x_{n-1}$.\\

\centerline{$St(\Omega):=\{x+ye_{n}:x+ze_{n}\in \Omega$ for some z and $|y|\leq\frac{1}{2} |\Omega\cap L_{x}| \}$.}

\subsubsection{Circular symmetrization}
The most popular method for symmetrization in the plane is $P\acute{o}lya$'s circular symmetrization. Let $\Omega \subset \mathbb{R}^{n+2}$ be a domain, then we denote the \textit{circularly symmetrized} $\Omega$ as $Cir(\Omega)$ and for $n\geq 3$ we denote the \textit{spherically symmetrized} $\Omega$ as $Sph(\Omega)$.\\
Let $\Omega \subset \mathbb{R}^{n+2}$ be a domain. For each $r\in (0,\infty)$ let $\Omega(r)=\{x\in \mathbb{S}^{n+1}:rx\in \Omega\}$. If $\Omega(r)=\mathbb{S}^{n+1}$, then the intersection of $Sph(\Omega)$ with the sphere $|x|=r$ is the full sphere, and if $\Omega(r)$ is empty then so is the intersection of $|x|r$. If $\Omega(r)$ is a proper subset of $\mathbb{S}^{n+1}$ and surface area $\sigma(\Omega(r))=A$ then the intersection of $Sph(\Omega)$ with $|x|=r$ is the cap $C(\theta):=\{(r,\phi):0\leq \phi\leq \theta\}$, where $\theta$ satisfies $\sigma(C(\theta))=A$. Moreover, $0\in Sph(\Omega)$ if and only if $0\in \Omega$.

\subsubsection{Polarization}
Let $\Omega\subset\mathbb{R}^{n}$ be a domain and $H^{n-1}\subset\mathbb{R}^{n}$ be a hyperplane through the origin. Denote the reflection across that plane as $\sigma_{H}$ or just $\sigma$ when it is clear from the context. Also, we denote the reflected $\Omega$ across hyperplane H as $\sigma \Omega$. Then, we denote the \textit{polarized} $\Omega$ as $\Omega^{\sigma}$ and define it as follows
\begin{equation*}
\Omega^{\sigma}\ni x^{\sigma}:=\left\{\begin{matrix}\sigma x & x\in (\Omega \setminus \sigma \Omega)\cap \mathbb{H}^{-}\\ x & x\in  (\Omega \cup \sigma \Omega)\cap \mathbb{H}^{+}\end{matrix}\right.
\end{equation*}

\begin{figure}[h]
\centering
\begin{minipage}[t]{.5\textwidth}
  \centering
  \includegraphics[scale=0.22]{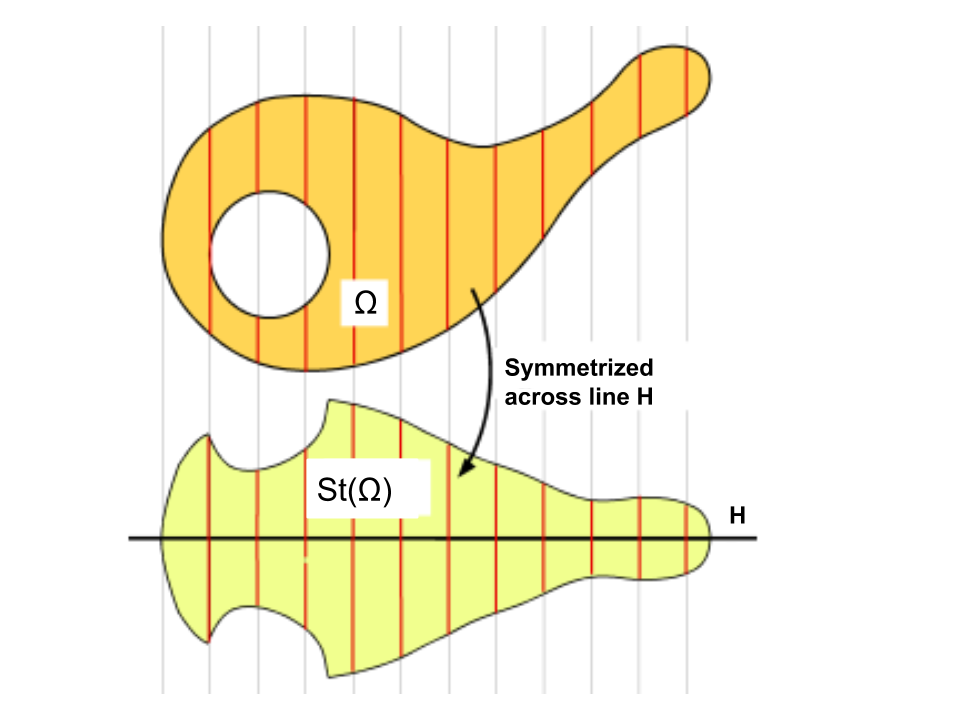}
  \caption{Steiner symmetrization of plane set $\Omega$}
  \label{fig:test1}
\end{minipage}%
\begin{minipage}[t]{.5\textwidth}
  \centering
  \includegraphics[scale=0.22]{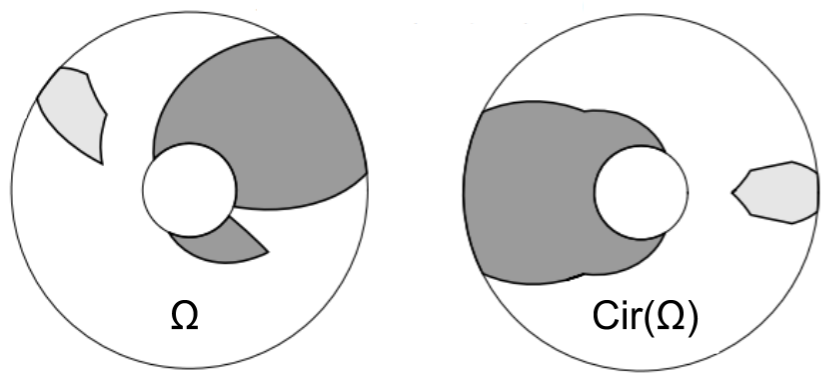}
  \caption{Circular symmetrization of plane set $\Omega$}
  \label{fig:test2}
\end{minipage}

\end{figure}

\begin{figure}[h]
\centering
\includegraphics[scale=0.3]{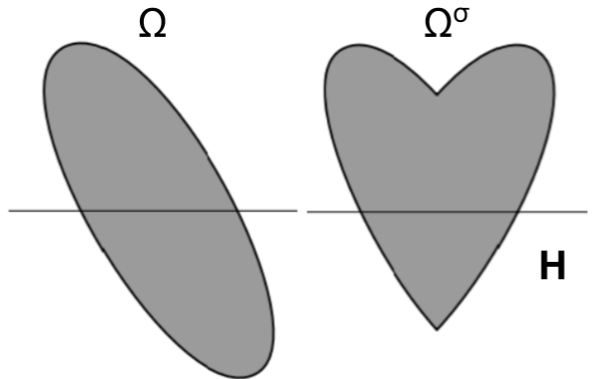}
\caption{Polarization of plane set $\Omega$}
\end{figure}

\pagebreak
\subsection{Brownian motion and Markov property }
The Markov property intuitively states that knowing the current position of a random process yields as much information as knowing the entire history of positions up to that point. A filtration on a probability space $(\Omega, \mathcal{F}, P)$ is a family $\{\mathcal{F}(t):t\geq 0\}$ of $\sigma-$ algebras such that $\mathcal{F}(s)\subset \mathcal{F}(t)\subset \mathcal{F}$ for $s\leq t$. Consider probability space, where $B_{t}$ is $\mathcal{F}_{t}$ measurable. Given the filtration $\mathcal{F}_{t}$ up to time t, then the Markov property for time $s$ says
\begin{equation*}
P(B_{t+s}\leq y|\mathcal{F}_{t})= P(B_{t+s}\leq y|B_{t})
\end{equation*}
or equivalently for starting point x 
\begin{equation*}
E_{x}(B_{t+s}|\mathcal{F}_{t})=E_{B_{t}}(B_{t+s}). 
\end{equation*}
For event A, we call $T_{A}:=\inf\limits_{t\geq 0}\{B_{t}\in A\}$ a \textit{stopping time}. Strong Markov property is similar to the Markov property, except that in the definition a fixed time t is replaced by a stopping time. The Strong Markov property for BM was proved by Hunt \cite{Hunt_1956} and for time s it says 
\begin{equation*}
E_{x}(B_{T_{A}+s}|\mathcal{F}_{T_{A}})=E_{B_{T_{A}}}(B_{T_{A}+s}).
\end{equation*}
Let $\Omega\subset\mathbb{R}^{n}$ be a closed or open set, then $\tau_{\Omega}:=\inf\limits_{t\geq 0}\{B_{t}\in \Omega\}$ will denote the \textit{hitting time} if $B_{0}=x\in \Omega^{c}$. On the other hand, for $B_{0}=x\in \Omega$, the $T_{\Omega}:=\inf\limits_{t\geq 0}\{B_{t}\in \partial \Omega\}$ will denote the \textit{exit time}, which we may also write as $T_{\partial \Omega}$. For more details see \cite{Peres_Morters_2010}.
\subsection{Harmonic measure and Exit probability}
In this section we will prove a connection of Brownian motion and complex analysis discovered by Kakutani (1944) i.e. equality of exit probability and harmonic measure. The harmonic measure appears in the Dirichlet problem: Let $\Omega\subset\mathbb{R}^{n}$ be a bounded domain and $\phi:\partial \Omega\to \mathbb{R}$ a continuous function, then the \textit{Dirichlet problem} is
\begin{equation*}
\Delta u=0 ~\text{and}~u|_{\partial \Omega}=\phi
\end{equation*}
Zaremba (1911) and Lebesgue (1924) gave examples of $\Omega$ where there is no solution. A sufficient condition is the $Poincar\acute{e}~ cone~ condition$: For each $x\in \partial \Omega$ there exists a cone $C_{x}(\alpha)$ based on it with opening angle $\alpha>0$ and for some $h>0$ it holds that $C_{x}(\alpha)\cap B_{h}(x)\subset \Omega^{c}$. All our sets in this thesis will satisfy this condition. Thus, if u is a solution, then for fixed $x\in \Omega$ Riesz representation theorem and the maximum principle yields a probability measure $\omega(x,\Omega)$ on $\partial \Omega$ s.t.
\begin{equation*}
u(x)=\int_{\partial \Omega}\phi(y)d\omega(x,\Omega)(y)
\end{equation*}
This measure $\omega(x,\Omega)$ is called the \textit{harmonic measure}. Thus, for any Borel subset $E\subset \partial \Omega$, the harmonic measure $\omega(x,\Omega)(E)$ is equal to the value at x of the solution to the Dirichlet problem with boundary data equal $\phi(y)=1_{E}(y)$ (characteristic functions can be approximated by continuous functions). Kakutani (1944) showed that 
\begin{thm}[Kakutani's theorem]
For notation as above
\begin{equation*}
\omega(x,\Omega)(E)=P_{x}(T_{E}=T_{\partial \Omega}).
\end{equation*}
\end{thm}
In other words, the probability that a Brownian motion starting at x will exit boundary $\partial \Omega$ via subset $E\subset\partial \Omega$ equals the harmonic measure of the Dirichlet problem $\Delta u=0$ and $u|_{\partial \Omega}=1_{E}$. The proof of this theorem can be found in \cite[section 3]{Peres_Morters_2010} and it is as follows: show that $u(x):=E_{x}[\phi(B_{T_{\Omega}})   ]=\int_{\partial \Omega} \phi(y) dP_{x}(B_{T_{\Omega} }=y)$ is harmonic in $\Omega$ and continuous in $\overline{\Omega}$.\\
Then the equality follows: For $x\in \partial \Omega$, $u(x)=E_{x}[\phi(B_{T_{\Omega}})  ]=E_{x}[\phi(B_{0})    ]=\phi(x)$, maximum principle and continuity yields 
\begin{equation*}
u(x)=\int_{\partial \Omega}\phi(y)d\omega(x,\Omega)(y)\Rightarrow \omega(x,\Omega)(\cdot)=P_{x}(B_{T_{\Omega} }=B_{T_{\cdot} }).
\end{equation*}
First, harmonicity of $E_{x}[\phi(B_{T_{\Omega}}) ]$. Consider ball $B_{\delta}(x)\subset \Omega$, then double conditioning and strong Markov property implies that
\begin{align*}
u(x)=&E_{x}[E_{x}[\phi(B_{T_{\Omega}}) |\mathcal{F}_{T_{B_{\delta}(x)}}] ]\\
=&E_{x}[u(B_{T_{B_{\delta}(x)}})]=\int_{\partial B_{\delta}(x)}u(y)d\sigma(y),
\end{align*}
where $\sigma$ is the uniform distribution on $\partial B_{\delta}(x)$. The last equality follows from the rotational invariance of the transition probability of Brownian motion. Thus, $u(x)$ satisfies the mean value property. Secondly, we will show continuity of $u(x)$ in $\partial \Omega$. This will follow from the Poincar$\acute{e}$ cone condition. \\
For any $z\in \partial \Omega$, by the Poincar$\acute{e}$ cone condition there is a cone $C_{z}(\alpha)$ and $h>0$ with $C_{z}(\alpha)\cap B_{h}(z)\subset \partial\Omega$. 
\begin{figure}[h]
\centering
\includegraphics[scale=0.3]{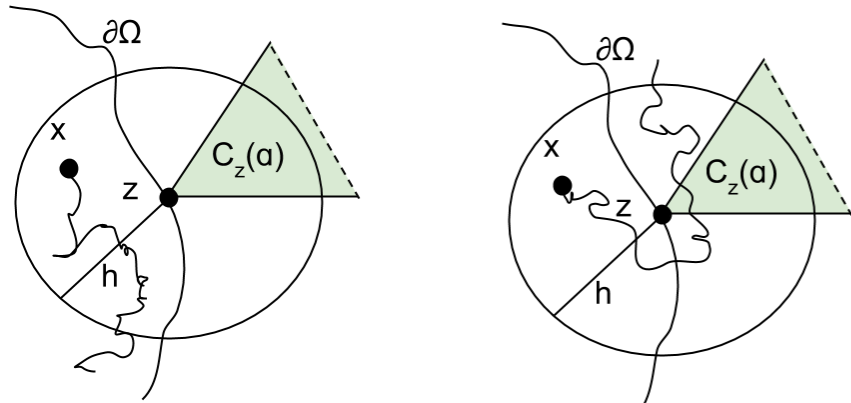}
\caption{Brownian motion avoiding the cone and hitting it before sphere $\partial B_{\delta}(z)$}
\label{split}
\end{figure}

Therefore, if the Brownian motion hits the cone $C_{z}(\alpha)$ before the sphere $\partial B_{\delta}(z)$, for any $\delta<h$ then $|z-B(T_{\partial \Omega})|<\delta$. By continuity of $\phi$ given $\varepsilon>0$, there is a $0<\delta\leq h$ s.t. $|\phi(y)-\phi(z)|<\varepsilon$ for all $y\in \partial \Omega$ with $|y-z|<\delta$. Therefore, $|z-B(T_{\partial \Omega})|<\delta\Rightarrow |\phi(B_{T_{\Omega} }) -\phi(z)|<\varepsilon$. By triangle inequality, for $z\in \partial \Omega$ and $x\in \overline{\Omega}$ with $|x-z|<2^{-k}\delta$ (for any integer k)
\begin{equation*}
|u(x)-u(z)|=|E_{x}\phi(B_{T_{\Omega} }) -\phi(z)|\leq E_{x}|\phi(B_{T_{\Omega} }) -\phi(z)|
\end{equation*}
then we split the last term into two events: the event that BM hits sphere $\partial B_{\delta}(z)$ before cone $C_{z}(\alpha)$ and the event that BM hits the cone $C_{z}(\alpha)$ before sphere $\partial B_{\delta}(z)$ i.e.
\begin{align*}
E_{x}|\phi(B_{T_{\Omega} }) -\phi(z)|=& E_{x}[|\phi(B_{T_{\Omega} }) -\phi(z)||T_{\partial B_{\delta}(z) }<T_{C_{z}(\alpha)} ]  P_{x}(T_{\partial B_{\delta}(z) }<T_{C_{z}(\alpha)})\\
+& E_{x}[|\phi(B_{T_{\Omega} }) -\phi(z)|T_{\partial \Omega }<T_{\partial B_{\delta}(z)}]    P_{x}(T_{\partial \Omega }<T_{\partial B_{\delta}(z)})\\
\leq & 2\left \| \phi \right \|_{\infty} P_{x}(T_{\partial B_{\delta}(z) }<T_{C_{z}(\alpha)})+\varepsilon P_{x}(T_{\partial \Omega }<T_{\partial B_{\delta}(z)}).
\end{align*}

Thus, it suffices to show that $P_{x}(T_{\partial B_{\delta}(z) }<T_{C_{z}(\alpha)})<\varepsilon$. 
\begin{lem}
Let $0<\alpha<2\pi$ and $C_{0}(\alpha)\subset\mathbb{R}^{n}$ be a cone based at the origin with opening angle $\alpha$, and $M=\sup \limits_{x\in \overline{B_{\frac{1}{2}}(0)}}P_{x}\{T_{\partial B_{1}(0)}< T_{C_{0}(\alpha)}   \}$, then $M<1$ and for any positive integer k and $h>0$ and $x,z\in \mathbb{R}^{n}$ s.t. $|x-z|<2^{-k}h$
\begin{equation*}
P_{x}(T_{\partial B_{h}(z)}<T_{C_{z}(\alpha)  })\leq M^{k}
\end{equation*}
\end{lem}
\begin{figure}[h]
\centering
\includegraphics[scale=0.3]{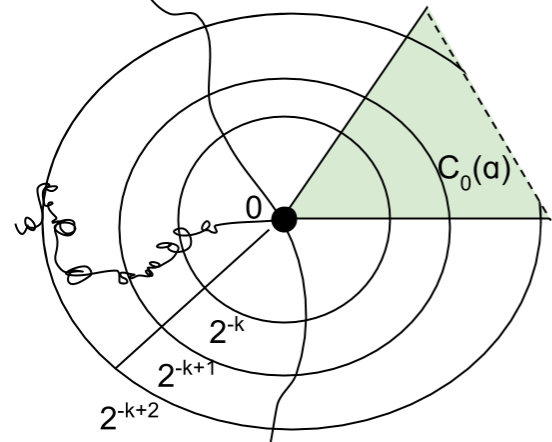}
\caption{Brownian motion avoiding the cone $C_{0}(\alpha)$}
\end{figure}

\begin{proof}
Brownian motion can be constructed s.t. $sup_{t\geq 0}B_{t}$ is arbitrarily close to $B_{1}(0)$. Then since $\alpha<2\pi$, we get $M<1$. If $x\in B_{2^{-k}}(0)$, then by the strong Markov property
\begin{equation*}
P_{x}(T_{\partial B_{1}(0)}<T_{C_{0}(\alpha)})\leq \prod_{i=0}^{k-1}\sup\limits_{x\in B_{2^{-k+i}}(0)  } P_{x}(T_{\partial B_{2^{-k+i+1}}(0)  }<T_{C_{0}(\alpha)})=M^{k}. 
\end{equation*}
Therefore, for any positive integer k and $h>0$, we have by scaling $P_{x}(T_{\partial B_{h}(z)}<T_{C_{z}(\alpha) })\leq M^{k}$, for all x with $|x-z|<2^{-k}h$.

\end{proof}

\pagebreak

\section{Brascamp, Lieb and Luttinger Inequality}

The Brascamp-Lieb-Luttinger inequality (BLL) provides a powerful and elegant method for obtaining and extending many of the classical geometric and physical isoperimetric inequalities of G. $P\acute{o}lya$ and $G. Szeg\ddot{o}$. Perhaps the most famous example of these inequalities is the celebrated Faber-Krahn inequality. That is, if $\lambda_{1}(D)$ and $\lambda_{1}(D^{*})$ are the first eigenvalues for the Laplacian with Dirichlet boundary conditions in $D$ and $D^{*}$, respectively, then
\begin{equation*}
\lambda_{D^{*}}\leq \lambda_{D}.
\end{equation*}
We will prove this in the "Principal eigenvalue" section using Kac's formula. Luttinger provided a new method, based on the Feynman-Kac representation of the heat kernel in terms of multiple integrals to prove the FK inequality\cite{Luttinger_1973}. The following inequality, proved by Brascamp, Lieb and Luttinger is a refinement of the original inequality of Luttinger \cite{Brascamp_Lieb_Luttinger_1974}. We also give a version on the sphere proved by Burchard and Schmuckenschl$\ddot{a}$ger \cite{Burchard_Schmuckenschlager_2001}. They also exist BLL inequalities in terms of inradius \cite{Banuelos_2001}. In this section, we will use BLL to prove the isoperimetric for exit probability and Riesz capacity for symmetric $\alpha-$stable processes of order $\alpha\in (0,2)$.

\begin{thm}\label{BLLinequalities}[Brascamp-Lieb-Luttinger inequalities]\\

\begin{enumerate}
\item\label{BLL}  
Let $\{f_{i}\}_{1\leq i\leq m}$  be nonnegative functions in $\mathbb{R}^{n}$ and $A\subset \mathbb{R}^{n}$ finite volume domain, then for any $z_{0}\in \mathbb{R}^{n}$\\
\centerline{$\int_{A^{m}}\prod_{i=1}^{m}f_{i}(z_{i}-z_{i-1})dz_{1}\cdots dz_{m}\leq \int_{(A^{*})^{m}}\prod_{i=1}^{m}f_{i}^{*}(z_{i}-z_{i-1})dz_{1}\cdots sz_{m}$.}
\item \label{BLL on sphere}[BLL on the sphere]\\
Let $\{A_{i}\}_{1\leq i\leq n}\subset \mathbb{S}^{n}$ be Borel sets and $\psi_{ij}:\mathbb{S}^{n}\times\mathbb{S}^{n}\to \mathbb{R}_{+}$ be non-increasing functions then \cite{Burchard_Schmuckenschlager_2001} 

\begin{align*}
\int_{(\mathbb{S}^{n})^{n}}\prod_{1\leq i\leq n}1_{x_{i}\in A_{i}}\prod_{1\leq j\leq n}\psi_{ij}(x_{i},x_{j})\prod_{1\leq i\leq n} d\sigma(x_{i})\\
\leq \int_{(\mathbb{S}^{n})^{n}}\prod_{1\leq i\leq n}1_{x_{i}\in A_{i}^{*}}\prod_{1\leq j\leq n}\psi_{ij}(x_{i},x_{j})\prod_{1\leq i\leq n} d\sigma(x_{i})
\end{align*}

\item \label{FL} [Friedberg-Luttinger inequality  \\
Let $\{F_{i}\}^{n}:\mathbb{R}^{n}\to [0,1]$ and $\{H_{i}\}^{n}$ nonegative nonincreasing radially symmetric functions in $\mathbb{R}^{n}$, then for $z_{m+1}:=z_{0}$ \cite{Friedberg_Luttinger_1976}]
\begin{align*}
\int_{\prod^{m}_{0} \mathbb{R}^{n}}[1-\prod_{i=0}^{m}(1-F_{i}(z_{i}))]\prod_{i=0}^{m}H_{i}(z_{i}-z_{i-1})dz_{0}\cdots dz_{m}\geq &\\ \int_{\prod^{m}_{0} \mathbb{R}^{n}}[1-\prod_{i=0}^{m}(1-St(F_{i})(z_{i}))]\prod_{i=0}^{m}H_{i}(z_{i}-z_{i-1})dz_{0}\cdots dz_{m}\geq &\\ \int_{\prod^{m}_{0} \mathbb{R}^{n}}[1-\prod_{i=0}^{m}(1-F_{i}^{*}(z_{i}))]\prod_{i=0}^{m}H_{i}(z_{i}-z_{i-1})dz_{0}\cdots dz_{m}
\end{align*}

\end{enumerate}

\end{thm}

\subsection{Symmetrization decreases exit probability}
Using the BLL inequalities we show that the exit probability of BM after time t ie. $P_{z}(T_{D}> t )$ increases with symmetrization. In other words, it becomes harder to escape when the domain is symmetrized. We will prove it for more general processes called n-dimensional symmetric $\alpha-$stable. Let $X_{t}$ be a n-dimensional symmetric $\alpha-$stable process of order $\alpha\in (0,2)$. Such a process has right continuous sample paths and stationary independent increments. Its infinitesimal generator is
\begin{equation*}
(-\Delta)^{\frac{\alpha}{2}}.
\end{equation*}
When $\alpha=2$ the process $X_{t}$ is just a n-dimensional Brownian motion $B_{t}$ running at twice the speed. Also, $X_{t}:=B_{2\sigma_{t}}$, where $\sigma_{t}$ is a stable subordinator of index $\frac{\alpha}{2}$ that is independent of $B_{t}$. Thus
\begin{equation*}
p_{\alpha}(t,x,y):=\int_{0}^{\infty}\frac{1}{(4\pi u)^{\frac{n}{2}}}e^{\frac{-|x-y|^{2}}{4u}}g_{\frac{\alpha}{2}}(t,u)du ,
\end{equation*}
where $g_{\frac{\alpha}{2}}(t,u)$ is the transition density of $\sigma_{t}$. Hence for every positive t, $p_{\alpha}(t,x,y)=f_{t}^{\alpha}(|x-y|)$ and the function is $f_{t}^{\alpha}(r)$ is decreasing. Thus, the conditions for BLL are satisfied. As with BM, we define as $P_{z}(T_{D,\alpha}> t )$ the exit probability of $X_{t}$. 
\begin{thm}
Let $D\subset \mathbb{R}^{n}$ be a domain of finite volume ,  $0<\alpha\leq 2$ , $z\in D$  and $t>0$ then\\
\centerline{$P_{z}(T_{D,\alpha}> t )\leq P_{0}(T_{D^{*},\alpha}> t )$.}
\end{thm}

\begin{proof}
By the right continuity of the sample paths and the Markov property of stable processes, we have
\begin{align*}
 P_{z } (T_{D,\alpha}\geq t)=& lim_{n\to \infty}P_{z} (X_{\frac{jt}{m}}\in D,j=1,...,m)\\
=& lim_{m\to \infty}\int_{D}\cdots \int_{D} \prod_{j=1}^{m} p_{\alpha}(\frac{t}{m},z_{j}-z_{j-1})dz_{1}\cdots dz_{m}\\
\leq& lim_{m\to \infty}\int_{D^{*}}\cdots \int_{D^{*}} p_{\alpha}(\frac{t}{m},0)\prod_{j=2}^{m} p_{\alpha}(\frac{jt}{m},z_{j}-z_{j-1})dz_{1}\cdots dz_{m}\\
=& lim_{m\to \infty}P_{0 } (X_{\frac{jt}{m}}\in D^{*},j=1,...,m)\\
=&P_{0 } (T_{D^{*}},\alpha\geq t).
\end{align*}

\end{proof}

Remark: Then it follows for closed sets by taking decreasing open sets. Similarly, it follows for $F_{\sigma}$ sets $D\subset \mathbb{R}^{n}$ (i.e. countable union of closed sets). Also, by the tail formulation $E_{x}(T_{D})=\int_{0}^{\infty}P_{x}(T_{D}>t)dt\leq E_{x}(T_{D^{*}})$.

\subsection{Symmetrization decreases $\alpha$-Riesz capacity }
The electrostatic capacity of an object is defined by the following problem. Assume the object is conducting and charged so that its surface has a constant (unit) potential, and the potential outside the object decays to zero at infinite distance. The capacity can then be defined in terms of the asymptotic decay at large distances of the solution to Laplace's equation in the space surrounding the object. For dimension 2, this is is called the logarithmic capacity and for $n\geq 3$ the Newtonian capacity (for details see \cite{Landkoff_1972}). As above we will prove the isoperimetric for more general capacities corresponding to n-dimensional symmetric $\alpha-$stable process $X_{t}$. The $\alpha$-Riesz kernel is
\begin{equation*}
k_{\alpha}(x-y)=\frac{\Gamma(n-\frac{\alpha}{2})}{\Gamma(\frac{\alpha}{2}) \pi^{\frac{n}{2} 2^{\alpha-1}  }}\frac{1}{|x-y|^{n-\alpha}} ,
\end{equation*}

where $n\geq 2$ and $0<\alpha<n$. Let A be a compact non-polar set in $\mathbb{R}^{n}$, the $\alpha-$Riesz capacity of A is defined by
\begin{equation*}
Cap_{\alpha}(A):=[\inf\limits_{\mu}\int\int k_{\alpha}(x-y) d\mu(x)d\mu(y)  ]^{-1}
\end{equation*}
where the infimum is taken over all probability Borel measures supported in A. If $\alpha=2$ and for $n\geq 3$,  this is the Newtonian capacity. In Getoor  \cite{Getoor_1965} it is proven that $Cap_{\alpha}(A)=lim_{t\to \infty}\frac{\int P_{z_{0}}(T_{A^{*},\alpha}\leq t)dz_{0}}{t}$. This is inspired by the same result for BM proved by Spitzer \cite{Spitzer_1964}(we prove this formula for BM in "Spitzer's formula" section). Further details for sharper asymptotic results such as $\int P_{z_{0}}(T_{A^{*},\alpha}\leq t)dz_{0}=Cap(A)t+\frac{1}{2\pi^{\frac{3}{2}}}Cap(A)^{2}+o(t^{\frac{1}{2}})$ can be found in \cite{Van_Berg_2007}.

\begin{thm}
Let $\alpha\in (0,2)$ and $A\subset \mathbb{R}^{n}$ be a bounded $F_{\sigma}$ set s.t. $vol(A)>0$ then 
\begin{equation*}
Cap_{\alpha}(A)\geq Cap_{\alpha}(St(A))\geq Cap_{\alpha}(A^{*}).
\end{equation*}
\end{thm}

\begin{proof}\cite{Mendez_2006}
Let $A_{k}$ be a decreasing sequence of compact sets such that the interior of $A_{k}$ contains A for all k and $\bigcap_{k=1}^{\infty}A_{k}=A$. Also, let $T_{A,\alpha}$ denote the stopping time for $X_{t}$. By the right continuity of the sample paths and the Markov property of stable processes, we have 
\begin{align*}
\int P_{z_{0} } (T_{A,\alpha}\leq t)dz_{0}=&\int 1-P_{z_{0} } (T_{A,\alpha}< t)dz_{0}=\int 1-P_{z_{0} } (X_{s}\in A^{c},0\leq s\leq t)dz_{0}\\
=&lim_{k\to \infty}lim_{n\to \infty}\int 1-P_{z_{0} } (X_{\frac{jt}{m}}\in A^{c}_{k},j=1,...,m)dz_{0}\\
=&lim_{k\to \infty}lim_{n\to \infty}\int\cdots \int [1-\prod_{j=1}^{m} I_{A^{c}_{k}}(z_{j}) ]\prod_{j=1}^{m}p_{\alpha}(\frac{t}{m},z_{j}-z_{j-1})dz_{0}\cdots dz_{m}\\
=&lim_{k\to \infty}lim_{n\to \infty}\int\cdots \int [1-\prod_{j=1}^{m} 1-I_{A_{k}}(z_{j}) ]\prod_{j=1}^{m}p_{\alpha}(\frac{t}{m},z_{j}-z_{j-1})dz_{0}\cdots dz_{m},
\end{align*}
where $I_{A_{k}}$ is the indicator function of $A_{k}$. Since $f_{t}^{\alpha}(x)$ is nonincreasing and radially symmetric, we can take $H_{m}=1$ and $F_{0}=0$ in the FL inequality to obtain

\begin{align*}
&\geq lim_{k\to \infty}lim_{n\to \infty}\int\cdots \int [1-\prod_{j=1}^{m} 1-I_{St(A_{k})}(z_{j}) ]\prod_{j=1}^{m}p_{\alpha}(\frac{t}{m},z_{j}-z_{j-1})dz_{0}\cdots dz_{m}\\
&=\int P_{z_{0}}(T_{A^{*},\alpha}\leq t)dz_{0}=:E_{St(A),\alpha}(t)\\
&\geq lim_{k\to \infty}lim_{n\to \infty}\int\cdots \int [1-\prod_{j=1}^{m} 1-I_{A_{k}^{*}}(z_{j}) ]\prod_{j=1}^{m}p_{\alpha}(\frac{t}{m},z_{j}-z_{j-1})dz_{0}\cdots dz_{m}\\
&=\int P_{z_{0}}(T_{A^{*},\alpha}\leq t)dz_{0}=:E_{A^{*},\alpha}(t).
\end{align*}

The $E_{A,\alpha}(t)$ is called the energy of A. Finally, from $Cap_{\alpha}(A)=lim_{t\to \infty}\frac{E_{A,\alpha}(t)}{t}$ \cite{Getoor_1965} it holds that\\
\begin{equation*}
Cap_{\alpha}(A)=lim_{t\to \infty}\frac{E_{A,\alpha}(t)}{t}\geq lim_{t\to \infty}\frac{E_{St(A),\alpha}(t)}{t}=Cap_{\alpha}(St(A))\geq lim_{t\to \infty}\frac{E_{A^{*},\alpha}(t)}{t}=Cap_{\alpha}(A^{*}).
\end{equation*}

\end{proof}

\subsection{Research Problems}
\begin{enumerate}
\item Let $K\subset \mathbb{R}^{n}$ for $n\geq 3$ be a compact set with finite positive volume, then for $\alpha\in (2,n)$
\begin{equation*}
Cap_{\alpha}(K)\geq Cap_{\alpha}(K^{*}).
\end{equation*}
\item Weakening the conditions on the integrands satisfying the BLL inequality. For example, in \cite{Burchard_Hajaiej_2006} it is shown that
\begin{thm}
Let $\{f_{i}\}_{1\leq i\leq m}$  be nonnegative functions in $\mathbb{R}^{n}$ vanishing at infinity, then for any $z_{0}\in \mathbb{R}^{n}$
\begin{align*}
&\int_{(\mathbb{R}^{n})^{k}}\prod_{i=1}^{m}F(f_{1}(\sum_{i=1}^{k}a_{1,i}x_{i}),..., f_{m}(\sum_{i=1}^{k}a_{m,i}x_{i})   )dz_{1}\cdots dz_{k}\\
&\leq \int_{(\mathbb{R}^{n})^{k}}\prod_{i=1}^{m}F(f_{1}^{*}(\sum_{i=1}^{k}a_{1,i}x_{i}),..., f_{m}^{*}(\sum_{i=1}^{k}a_{m,i}x_{i})   )dz_{1}\cdots dz_{k}
\end{align*}
where F is the distribution function of a Borel measure $\mu$ on $\mathbb{R}^{n}$ i.e. $F(y_{1},...,y_{m})=\mu([0,y_{1}),...,[0,y_{m}))$ and $a_{ij}\in \mathbb{R}$. Such F are left-continuous, nonnegative and $\Delta_{i_{1},...,i_{l}}F\geq 0$ for all choices $\{i_{j}\}^{l}\subset (1,..,m)$.

\end{thm}

\end{enumerate}

\pagebreak
\section{Baernstein star-function}

The star-function method is used for the solution of certain extremal problems (i.e. minimizers/maximizers) for which the competing functions u are subharmonic and the
expected extremal function v is harmonic in a symmetric region. For each u, Baernstein defined a certain maximal
function $u^{\star}$, the star-function of u,\cite{Baernstein_1974} and showed that the solution of the extremal problem is reduced to the inequality $u^{\star} \leq v^{\star}$. The heart of the method is the fact that $v^{\star}$ remains subharmonic, while the symmetry of the extremal domain implies that $v^{\star}$ is harmonic. It follows that the function $u^{\star}-v^{\star}$ is subharmonic and therefore, in order to prove the desired inequality $u^{\star}-v^{\star}\leq 0$ one can use the maximum principle. In this section, we will exemplify this technique to prove the isoperimetric for the harmonic measure. \\

For $0<a<\infty$, let $g\in L^{1}([-a,a],\mathbb{R})$. Define the \textit{Baernstein $\star$-function} $g^{\star}:[0,a]\to \mathbb{R}$ by \cite{Baernstein_2002} \\
\centerline{$g^{\star}(l):=\sup\limits_{E\subset [-a,a]}\int_{E}g(s)ds$,}\\

where the supremum is taken over all Lebesgue measurable sets $E\subset[-a,a]$ with $|E|=2l$. For each $l$ by continuity of $\lambda(t):=|\{x:g(x)>t\}|$, one can show that the supremum is attained at set E and $g^{\star}(l)=\int_{-l}^{l}g^{*}(s)ds$, where $g^{*}$ be the decreasing symmetric rearrangement of g. So intuitively the star function measures the symmetric mass centered at the origin. For example, on the punctured disk $\mathbb{D}\setminus \{0\}$ the harmonic function $f(re^{i\theta})=log(r)$, has $f^{*}(re^{i\theta})=log(1-r^{2})$ and thus $f^{\star}(re^{i\theta})=log(1-r^{2})2\theta$.

\begin{figure}[h]
\centering
\includegraphics[scale=0.3]{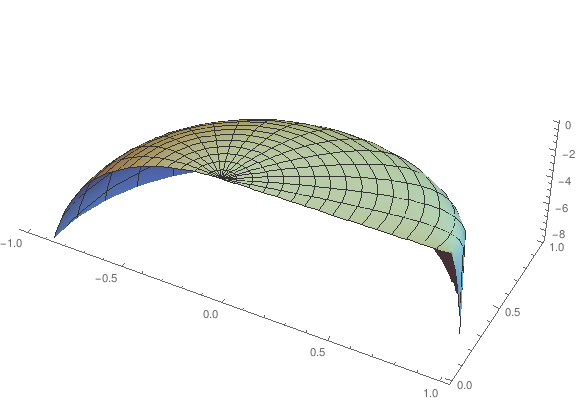}
\caption{Baernstein star function of $f(re^{i\theta})=log(r)$ on $\mathbb{D}\setminus \{0\}$}

\end{figure}

\subsection{Properties of Baernstein's star-function}
The following propositions are proved in \cite{Hayman_1989}.

\begin{prop}\label{symmetrization}[Baerstein star-function properties]
\begin{enumerate}
\item \label{symmetric} Let $g^{*}$ be the decreasing symmetric rearrangement of g, then $g^{\star}(\theta)=\int_{-\theta}^{\theta}g^{*}(s)ds$.

\item \label{equivalent} For $g,h\in L^{1}[-a,a]$, the following are equivalent:\\

a)For every convex increasing function $\Phi:\mathbb{R}\to \mathbb{R}$ holds
\begin{equation*}
\int_{-a}^{a}\Phi(g(s))ds\leq \int_{-a}^{a}\Phi(h(s))ds.
\end{equation*}

b)Let $g^{*},h^{*}$ be the decreasing symmetric rearrangement of g,h then for $s\in [0,a]$ holds
\begin{equation*}
g^{*}(s)\leq h^{*}(s).
\end{equation*}

\item If $g,h\in L^{1}[-a,a]$ and $g^{*}(s)\leq h^{*}(s)$ $\forall s\in [0,a]$, then 
\begin{equation*}
ess~sup_{[-a,a]}g\leq ess~sup_{[-a,a]}h.
\end{equation*}

\item \label{Baernstein's Fundamental Theorem} [Subharmonicity properties of star-function]\\
Suppose u is subharmonic in annulus $A_{r_{1},r_{2}}:=\{r_{1}<|z|<r_{2}\}\subset \mathbb{C}$ for $0\leq r_{1}<r_{2}\leq \infty$, then $u^{\star}$ is subharmonic in $A_{r_{1},r_{2}}\cap \mathbb{H}^{+}=\{r_{1}<|z|<r_{2}:0<argz<\pi\}$.\\

\item 
Suppose $u=u_{1}-u_{2}$, where $u_{i}$ are subharmonic in the disk $D_{R}:=\{|z|<R\}\subset \mathbb{C}$ and finite at the origin, then $u^{\star}(z)+\int_{-\pi}^{\pi}u_{2}(re^{i\theta})d\theta$ is subharmonic in $D_{R}\cap \mathbb{H}^{+}=\{|z|<R:0<argz<\pi\}$ and continuous in $\overline{D_{R}\cap \mathbb{H}^{+}}\setminus \{0\}$.

\end{enumerate}
\end{prop}

\subsection{Circular Symmetrization and Exit probability}

Suppose that D is a domain lying in $D_{R}:=\{|z|<R\}$ and let $\alpha$ be the intersection of the boundary of D with $|z|=R$. Let $Cir(D)$ be the circularly symmetrized domain of D and let $\alpha^{*}$ be the intersection of the boundary of $Cir(D)$ with $|z|=R$. We set
\begin{equation*}
u(z):=\omega(z,\alpha,D) and v(z):=\omega(z,\alpha^{*},Cir(D)).
\end{equation*}
we define $u=0$ outside D and $v=0$ outside $Cir(D)$. Then \cite{Hayman_1989}

\begin{thm}\label{mean_isoperimetric}[Isoperimetric of exit probability]\\
Let $\Phi:\mathbb{R}\to \mathbb{R}$ be convex non-decreasing function $\Phi$ and $r\in (0,1)$ , then 
\begin{equation*}
\int_{-\pi}^{\pi}\Phi(u(re^{i\theta}))d\theta \leq\int_{-\pi}^{\pi}\Phi(v(re^{i\theta}))d\theta.
\end{equation*}

then by the last proposition
\begin{equation*}
\sup\limits_{|z|=r}\omega(z,\alpha, D) \leq \sup\limits_{|z|=r}\omega(z,\alpha^{*},Cir(D)) =\omega(r,\alpha^{*}, Cir(D)).
\end{equation*}

In probabilistic terms, 
\begin{equation*}
\sup\limits_{|z|=r}P_{z}(T_{\alpha}=T_{\partial D})\leq P_{r}(T_{\alpha^{*}}=T_{\partial Cir(D)}) .
\end{equation*}
\end{thm}
\begin{figure}[h]
\centering
\includegraphics[scale=0.3]{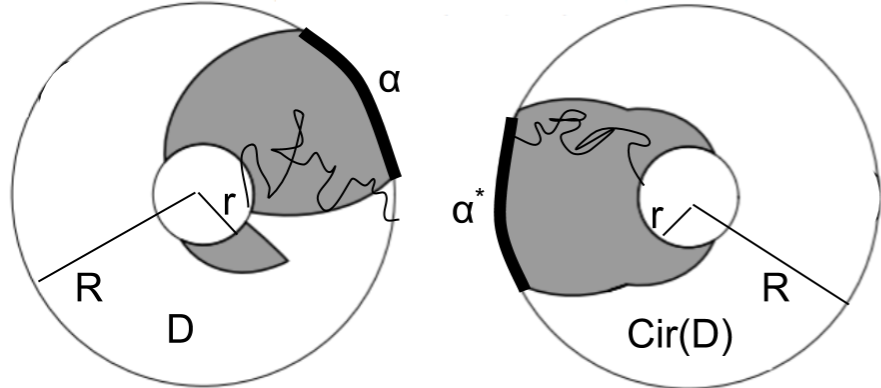}
\caption{Comparing the exit probability of BM starting from circle $|z|=r$}
\end{figure}

\begin{proof}
We sketch the proof by \cite{Hayman_1989} and elaborate the parts where the star function's properties are used. We assume the domain of D is smooth. For the general domain we use an expanding sequence of smooth domains $D_{n}$ whose union is D. For now we also assume that $\alpha=\{|z|=R\}$ i.e. the inner boundary $\beta:=\partial D\cap \{|z|<R\}$ does not intersect $\{|z|=R\}$. The harmonic measures u,v are subharmonic and continuous on $D_{R}$ . Also u and v are equal to 1 on $|z|=R$ and equal to 0 at points of $D_{R}$ outside D and $Cir(D)$ respectively. \\

\begin{figure}[h]
\centering
\includegraphics[scale=0.3]{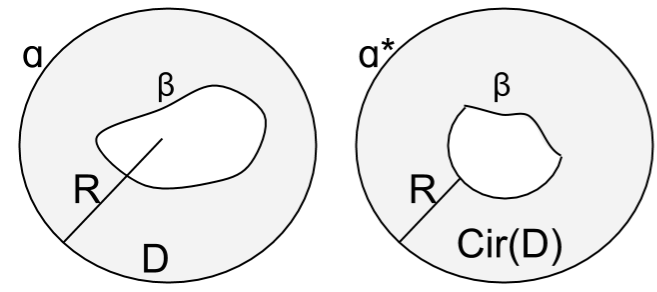}
\caption{Assuming that $\alpha=\{|z|=R\}$}
\end{figure}

From \pref{symmetrization}.\ref{Baernstein's Fundamental Theorem}, we have that since u is subharmonic in $D_{R}$, then $u^{\star}$ is subharmonic and continuous in $D_{R}\cap \mathbb{H}^{+}$ and thus in $Cir(D)\cap \mathbb{H}^{+}$. Similarly, using harmonicity of v and maximum principle, one can show that the  $v^{\star}$ is harmonic and continuous in $Cir(D)\cap \mathbb{H}^{+}$ (\cite[lemma 9.2.4]{ Hayman_1989}).\\
Thus, for positive $\varepsilon$ and $z=re^{i\theta}$ we have that $Q(z):=u^{\star}(z)-v^{\star}(z)-\varepsilon\theta$ is subharmonic in $D^{*+}:=Cir(D)\cap \mathbb{H}^{+}$ and continuous on $\overline{D_{R}\cap \mathbb{H}^{+}}$. Next we will show that for $z\in D_{R}\cap \mathbb{H}^{+}$ 
\begin{equation*}
Q(z)=u^{\star}-v^{\star}-\varepsilon\theta\leq 0
\end{equation*}

and so letting $\varepsilon\to 0$ we deduce the claim of the theorem due to the equivalence from \pref{symmetrization}.\ref{equivalent}.\\

We suppose that M is the supremum of Q(z) in $D^{*+}$. We will assume that $M>0$ and obtain a contradiction. Since Q(z) is subharmonic in $D^{*+}$, by the Maximum principle there exists a point $\zeta$ on the boundary of $D^{*+}$ and a sequence $z_{n}\in D^{*+}$ s.t. \begin{equation*}
z_{n}\to \zeta, Q(z_{n})\to M as n\to \infty.
\end{equation*}
Next we split into cases on the point $\zeta$ where the maximum is attained. When $\zeta$ is on the circle $|z|=R$ (case 1), on the positive axis (case 2), along the circle $|z|=r<R$ (case 3 and 4), $\zeta=0$ (case 5) and $\zeta=\infty$ (case 6).\\

Case 1: Assume $|\zeta|=R$. On $|z|=R$, then by definition of harmonic measures $u(z)=v(z)=1$ and so $u^{\star}(Re^{i\theta})=v^{\star}Re^{i\theta}=2\theta$. Thus, $Q(Re^{i\theta})=-\varepsilon\theta<0$ and so this contradicts $M>0$. So assume $|\zeta|<R$ for the rest of the cases.\\

Case 2: On the positive real axis by definition $u^{\star}(z)=u^{\star}(re^{i\cdot 0})=0$; so $M=Q(\zeta)=0$ and in turn $\zeta$ cannot lie there.\\

For case 3 and 4 we suppose that $\zeta=re^{i\phi}$, where $0<\phi\leq \pi$.\\

Case 3: If the circle $|z|=r$, does not intersect D, then $u^{\star}(\zeta)=v^{\star}(\zeta)=0$ and so $M<0$.\\

Case 4: If the circle $|z|=r$ intersects D in a set of measure 2m, then $2m<2r\phi$ and u(z)=v(z)=0 at points of $|z|=r$ outside $D\cap \{|z|=r\}$. By \pref{symmetrization}\ref{symmetric},
\begin{equation*}
u^{\star}(re^{i\phi})=\int_{-\phi}^{\phi}u^{*}(re^{i\psi}) d\psi,
\end{equation*}
and so by the continuity of the decreasing rearrangements,
\begin{align*}
\frac{\partial}{\partial \psi}u^{\star}(re^{i\psi})|_{\psi=\phi}&=2u^{*}(\phi)=0  \\
\frac{\partial}{\partial \psi}v^{\star}(re^{i\psi})|_{\psi=\phi} &=2v^{*}(\phi)=0.
\end{align*}
Thus,
\begin{equation*}
\frac{\partial}{\partial \psi}Q(re^{i\psi})|_{\psi= \phi} =-\varepsilon<0,
\end{equation*}
and so Q is decreasing i.e. for small positive h
\begin{equation*}
Q(re^{i(\phi-h)})>Q(re^{i\phi})\geq M.
\end{equation*}

This contradicts our assumption that M is the supremum of Q in $D^{*+}$. The same argument applies if $\phi=\pi$ and not the whole circle $|z|=r$ in D, so that this circle meets the complement of D in a non-empty set of measure zero.\\

Case 4: We suppose next that $\zeta=-r$, and that the whole circle $|z|=r$ lies in D and so in $Cir(D)$. Let $\rho_{1}<|z|<\rho_{2}$ be the largest annulus that is contained in D and contains $|z|=r$. Then u, v are harmonic in  $\rho_{1}<|z|<\rho_{2}$, hence
\begin{equation*}
Q(-\rho)=\int_{-\pi}^{\pi}(u(\rho e^{i\theta})-v(\rho e^{i\theta}))d\theta-\varepsilon \pi
\end{equation*}
is a linear function of $log\rho$ for $\rho_{1}<\rho<\rho_{2}$ (by Green's theorem). Since $Q\leq M$ everywhere and attains its maximum $Q(\zeta)=Q(-r)=M$ inside the annulus, we deduce by maximum principle for harmonic $log(\rho)$ and continuity of $Q(-\rho)$ in $0\leq \rho\leq \infty$ that \\
\begin{equation*}
Q(-\rho)=M ~(\rho_{1}\leq \rho\leq \rho_{2}),
\end{equation*}
Also, D cannot consist of the punctured plane, so that either $0<\rho_{1}<\infty$ or $0<\rho_{2}<\infty$. Since D is open, it does not contain the whole circle $|z|=\rho_{1}$ or $|z|=\rho_{2}$, and so we are reduced to case 3 and again obtain a contradiction.\\

Case 5: This leaves the possibilities $\zeta=0$ or $\zeta=\infty$. If $\zeta=0$ is a boundary point of D then $u(0)=v(0)=0$ and so
\begin{equation*}
M=\lim\sup\limits_{z\to 0}Q(z)\leq 0,
\end{equation*}
contrary to hypothesis. We suppose then that $\zeta=0$ is an interior point of D and let $|z|<\rho$ be the largest disk, with center 0, that lies in D. Since $u(z)-v(z)$ are harmonic at $z=0\in D$, we deduce that, uniformly in $\theta$ as $r\to 0$,
\begin{equation*}
Q(re^{i\theta})=(2u(0)-2v(0)-\varepsilon)\theta+o(1)=:A\theta+o(1).
\end{equation*}

If $A\leq 0$, we deduce that\\
\begin{equation*}
0<M=\lim\sup\limits_{z\to 0}Q(z)\leq 0,
\end{equation*}
which gives a contradiction. If $A>0$,
\begin{equation*}
Q(re^{i\theta})\leq Q(-r)
\end{equation*}

for small r. So by maximum principle $M=Q(-r)=A\pi$ for $0<r\leq \rho$. So because $Q(-\rho)=M$ and $|z|=\rho$ is not contained entirely in D, this leads to a contradiction as in case 3.\\

Case 6: If $\zeta=\infty$, we argue similarly. If $\infty\in D$, then u,v are harmonic in a neighbourhood $|z|>\rho$ of $\infty$ and thus as in case 5 obtain $M=0$. If $\infty \notin D$, since $u=v=0$ in $D^{c}$, the limit $\lim\sup\limits_{z\to \infty}Q(z)\leq 0$, which is again a contradiction. Thus the existence of $\zeta$ always leads to a contradiction. Thus, we have the theorem if D is smooth. In the case $\alpha\subset\{|z|=R\}$, by considering domains $D_{n}=D\cup \{R_{n}<|z|>R\}$ one can obtain the result.\\

\begin{figure}[h]
\centering
\includegraphics[scale=0.3]{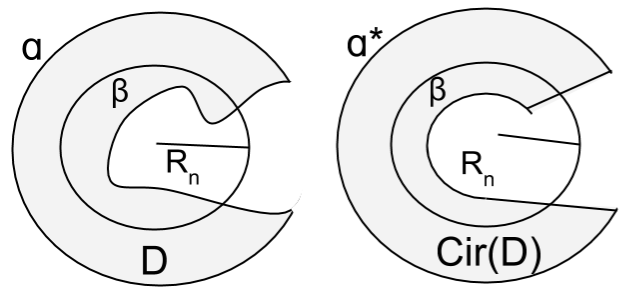}
\caption{Assuming that $\alpha\subset\{|z|=R\}$}
\end{figure}

\end{proof}
\pagebreak

\subsection{Steiner Symmetrization and Higher dimensions}

The idea is to transfer the analoguous problem in which integrals over sets on concetric circles are replaced by integrals over sets on parallel lines. 
For fixed $-\infty\leq x_{1}<x_{2}<\infty$, set\\
\begin{equation*}
B(x_{1},x_{2}):=\{x+iy\in \mathbb{C}: x_{1}<x<x_{2}\}~\text{ and}~ B^{+}(x_{1},x_{2}):=B(x_{1},x_{2})\cap \mathbb{H}^{+},
\end{equation*}
and for $u:B(x_{1},x_{2})\to\mathbb{R}$, define the \textit{vertical $\star$-function} $u^{\star}:B^{+}(x_{1},x_{2})$ by\\
\begin{equation}\label{definition}
u^{\star}(x+iy):=\sup\limits_{E}\int_{E}u(x+it)dt
\end{equation}

where the sup is taken over all $E\subset\mathbb{R}$ with $|E|=2y$. 

\begin{thm}
Let D be a domain with $D\subset B(-\infty,x_{2})$, where $x_{2}<\infty$. Set $L(x_{2}):=\{x+iy\in \partial D:x=x_{1} \}$. Then, for each convex increasing function $\Phi:\mathbb{R}\to\mathbb{R}$ and each $x\in (-\infty,x_{2})$ holds for Steiner symmetrization with respect to real axis
\begin{equation*}
\int_{\mathbb{R}}\Phi(\omega(x+iy,L(x_{2}),D))dy\leq\int_{\mathbb{R}}\Phi(\omega(x+iy,St(L(x_{2})),St(D))dy.
\end{equation*}

This implies $sup_{y}\omega(x+iy,L(x_{2}),D)\leq sup_{y} \omega(x+iy,St(L(x_{2})),St(D))$ or in probabilistic terms\\
\begin{equation*}
sup_{y}P_{x+iy}(T_{L(x_{2})}=T_{\partial D})  \leq sup_{y}P_{x+iy}(T_{St(L(x_{2})}=T_{\partial St(D)}).
\end{equation*}

\end{thm}

Remark: This result is also proved using Ahflor's distortion theorem and Brownian motion in Halliste \cite{Haliste_1965}.\\
The above results can be extended to higher dimensions \cite{Baernstein_Taylor_1976}. Much of the theory works like it does for Steiner and circular symmetrization in the plane. One significant change, though, is that if u is subharmonic with respect to the Laplace operator, the $u^{\star}$ will be subharmonic with respect to a possibly different operator, which depends on the symmetrization process. 

\subsection{Research Problems}
 Is the \thmref{mean_isoperimetric} true for multiply connected $Cir(D)$? For details see \cite{Baernstein_2002}.

\pagebreak
\section{Dubinin's Desymmetrization}
By contrast to the known symmetrization transformation, desymmetrization enables us to obtain ‘reverse’ estimates. Originally, desymmetrization was designed to solve Gonchar’s problem of harmonic measure \cite{Dubinin_1985}. Later it became clear that this transformation is interesting on its own  \cite{  Dubinin_2014}.

In this section, we will exemplify the properties of desymmetrization by proving Gonchar's problem. We will follow the second proof. Before, we state the problem, we need some notations. There are two proofs available; by Dubinin in \cite{Dubinin_2014} and by Baernstein \cite{Baernstein_1987}(nice exposition in \cite{Reid_1996}). \\
We start by defining a slit domain in $\mathbb{C}$. Let $K=[a,1]\subset [0,1]$, then a \textit{radial slit} is $zK:=\{\lambda z:\lambda\in K\}$. For $K=[a,1]$ and $\alpha:=(\alpha_{1},...,\alpha_{n})$, where $0\leq \alpha_{1}\leq \cdots \leq \alpha_{n}\leq 2\pi$, we call $\Omega_{\alpha}\subset\mathbb{D}$ a \textit{radially slit disk} if 
\begin{equation*}
\Omega_{\alpha}=\mathbb{D}\setminus \bigcup_{k=1}^{n}e^{i\alpha_{k}}K
\end{equation*}

In case, $\widetilde{\alpha}_{j}=\frac{j-1}{n}\pi$ (equally spaced angles), we will denote the radial slit disk as $\widetilde{\Omega}$ or $\Omega_{\widetilde{\alpha}}$.\\
Suppose, $\Omega_{\alpha}$ and $\Omega_{\widetilde{\alpha}}$ are two radially slit domains both with n slits formed by $\alpha=(\alpha_{1},...,\alpha_{n})$ and $\widetilde{\alpha}=(0,...,\frac{j-1}{n}\pi,...,\frac{n-1}{n}\pi)$ and  $K=[a,1]$. Let $S_{\alpha}=\bigcup_{k=1}^{n}e^{i\alpha_{k}}K$. For these we consider the harmonic measures starting from zero $\omega(0,S_{\alpha},\Omega_{\alpha})$ and $\omega(0,S_{\widetilde{\alpha}},\Omega_{\widetilde{\alpha}})$. Then we have the following theorem

\begin{figure}[h]
\centering
\includegraphics[scale=0.2]{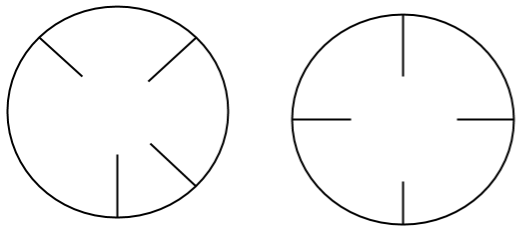}
\caption{$\Omega_{\alpha}$ and $\Omega_{\widetilde{\alpha}}$ slit domains respectively}
\end{figure}

\begin{thm}\label{Gonchar-Dubinin theorem}[Gonchar-Dubinin theorem]\\
Given $K=[a,1]\subset [0,1]$ and $\alpha:=(\alpha_{1},...,\alpha_{n})$ where $0\leq \alpha_{1}\leq \cdots \leq \alpha_{n}\leq 2\pi$, then for notation as above
\begin{equation*}\omega(0,S_{\alpha},\Omega_{\alpha})\leq  \omega(0,S_{\widetilde{\alpha}},\Omega_{\widetilde{\alpha}}),  \end{equation*}
or in probability terms
\begin{equation*}
P_{0}(T_{S_{\alpha}}=T_{\partial \Omega_{\alpha}})\leq P_{0}(T_{S_{\widetilde{\alpha}}}=T_{\partial \Omega_{\widetilde{\alpha}}}).
\end{equation*}

\end{thm}

Remark: The probability of exiting via one of the n- slits instead of the boundary of unit disk, increases when the n-slits are evenly placed across the disk.

\subsection{Background material}

In this section we describe the desymmetrization process. We will only mention Dubinin's version \cite{Dubinin_2014} to convey the main idea. Because we will use Baernstein's version in the proof, we will also mention its properties.

\subsubsection{Dubinin's definition}
Let $L^{*}_{k}$ , where $k=1,...,n$ and $n\geq 2$, be rays emanating from the origin at equal angles (eg. $\widetilde{\alpha}$). Let $\Phi$ be the group of symmetries of $\overline{\mathbb{C}}$ formed by the composites of the reflections in the rays $L^{*}_{k}$ and in the bisectors of the angles formed by these rays. Thus, we say a set $D\subset \mathbb{C}$ is \textit{$\Phi$-symmetric} if $\forall \phi\in \Phi$, $\phi(D)=D$. A function u on $\Phi-$symmetric D is \textit{$\Phi$-symmetric} if  $\forall \phi\in \Phi$ , $u(\phi(z))=u(z)$.\\
Let $A\subset \mathbb{D}$ and $u$ be $\Phi-$symmetric defined on $\Phi-$symmetric $\Omega\subset \mathbb{C}$ , we construct $A^{dis}$  and $u^{dis}$ as follows. 
\begin{enumerate}
\item We partition $\mathbb{C}$ into sets $\{P_{k}\}_{k=1}^{m}$ that are $\Phi-$invariant i.e. $\{\phi(P_{k})\}_{k=1}^{m}=\{P_{k}\}_{k=1}^{m}$ for all $\phi$. 
\item We consider angles $\{\theta_{k}\}_{k=1}^{m}$ s.t. the rotated $P_{k}$ , $S_{k}:=e^{i\theta_{k}}P_{k}$ satisfy $\mathbb{C}=\bigcup_{1}^{m}S_{k}$ and for each $S_{k}\cap S_{m}$ there exists an isometry $\phi\in \Phi$ s.t. $\phi(e^{-i\theta_{k}}(S_{k}\cap S_{m}))=e^{-i\theta_{m}}(S_{k}\cap S_{m})$.\\
Then we define 
\begin{equation*}
A^{dis}:=\bigcup_{1}^{m}e^{i\theta_{k}}(P_{k}\cap A) and u^{dis}(z)=u(e^{-i\theta_{k}}z) for z\in S_{k}\cap \Omega^{dis}.   
\end{equation*}

\end{enumerate}

\begin{figure}[h]
\centering
\includegraphics[scale=0.3]{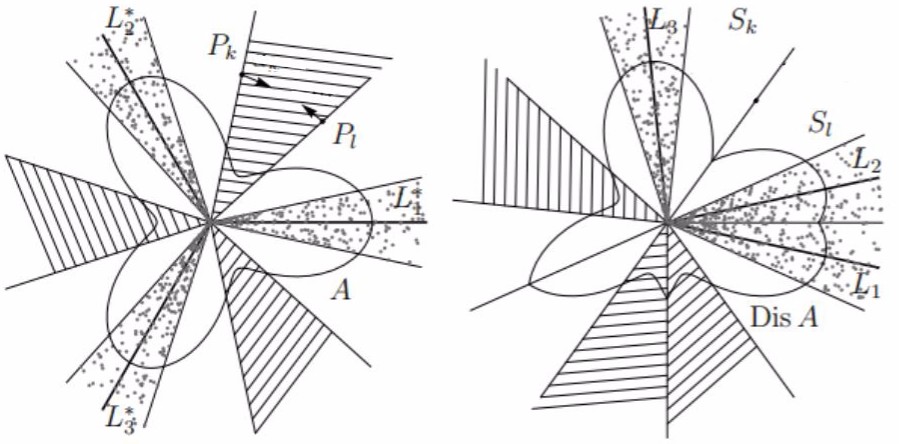}
\caption{Desymmetrization of domain A}
\end{figure}

In Gonchar's problem, we start with $\Omega_{\widetilde{\alpha}}$ (each of the equally spaced slits corresponding to a $L^{*}_{k}$) and then desymmetrize according to some angles $\theta_{k}$ such that the new slits will have angles $\alpha=(\alpha_{1},...,\alpha_{n})$.

\subsubsection{Desymmetrization properties}

Denote by $I_{\Omega}(u)$ the Dirichlet integral $I_{\Omega}(u):=\int_{D}|\bigtriangledown u(z)|^{2}dz$. A function $f:\mathbb{C}\to \mathbb{R}$ is \textit{admissible} if it is real piecewise smooth and Lipschitz. 

\begin{prop}\label{desy prop}[Desymmetrization properties]\\

The following properties are true for both Dubinin's and Baernstein's $f^{dis}$:
\begin{enumerate} 
\item Let f be an admissible function on $A(r_{1},r_{2})$ then $f^{dis}$ is also admissible on $A(r_{1},r_{2})$

\item If $\partial_{\theta}f(z)$ along each circle exists, then $\partial_{\theta}f^{dis}=\partial_{\theta}f(z)$ and so
\begin{equation*}
I_{A(r_{1},r_{2})}(f)=I_{A(r_{1},r_{2})}(f^{dis}).
\end{equation*}

The following properties are true for Baernstein's $f^{dis}$.\\

\item \label{equal on slits}$f^{dis}(re^{i\alpha_{j}})=f(re^{i\widetilde{\alpha_{j}} })$

\item f and $f^{dis}$ are equidistributed: $|\{f>t\}|=|\{f^{dis}>t \}|~\forall t$.

\item f and $f^{dis}$ have the same valence on each circle: if the equation $f(x)=y$ has m solutions, then so does $f^{dis}(x)=y$.%

\end{enumerate}
\end{prop}

\subsubsection{Dirichlet integral properties}

The following are some general properties about Dirichlet integrals and Harmonic measures we will need for the proof. Given annulus $A(r_{1},r_{2})\subset\mathbb{D}$, we denote by $D_{r_{1},r_{2}}(\theta_{1},\theta_{2})$ its \textit{sectors} $\{r e^{i\theta}\in A(r_{1},r_{2}): r_{1}\leq r\leq r_{2}, \theta_{1}\leq \theta\leq \theta_{2}\}$.\\

\begin{thm}\label{Dirichlet integral}[Properties of Dirichlet integral]\\

Let $\Omega\subset \mathbb{C}$ be a domain, $f,g$ admissible functions on $\Omega$.
\begin{enumerate}
\item \label{Modified Polya inequality}
Modified $P\acute{o}lya-Szeg\ddot{o}$ inequality\\
Let $f:\Omega\to [0,1]$ be admissible and on each circle $\{|z|=r\}$ equal to 0 and 1 only at n distinct points (i.e. $|f=0|=n=|f=1|$). Then, 
\begin{equation*} I_{\Omega}(f^{*}(re^{in\cdot \theta}))\leq I_{\Omega}(f), \end{equation*}
where $f^{*}$ is the \textit{circular symmetrization} of f and $f^{*}(re^{in\cdot \theta})$ is called the \textit{n-fold symmetrization}.\\

\item \label{harmonic dirichlet integral}
Let f,g be admissible on $\overline{\Omega}$,  $f=g $ on $\partial \Omega$ and f is harmonic on $\Omega$ then\\
\begin{equation*}
I_{\Omega}(f)\leq I_{\Omega}(g).
\end{equation*}

\item \label{Conformal invariance}
Conformal invariance of Dirichlet integrals\\
Let $\Phi:\Omega'\to \Omega$ be a conformal map and f admissible on $\Omega$, then $f\circ \Phi$ is admissible and 
\begin{equation*} I_{\Omega}(f)= I_{\Omega'}(f\circ \Phi) \end{equation*}

\item \label{annulus integral}
Dirichlet integrals of harmonic functions over annulus\\

Let $X:=\{g\in C^{1}(A(r_{1},r_{2})):\Delta g=0,~ g(r_{1}e^{i\theta})=1$ and $g(r_{2}e^{i\theta})=0$ $\forall \theta\in (0,2\pi) \}$, for some constants $r_{1}<r_{2}\leq 1$ i.e. functions that solve the Dirichlet problem with the above conditions. Then
\begin{equation*}I_{A(r_{1},r_{2})}(g)=2\pi (log(\frac{r_{2}}{r_{1}})^{-1}  )  \end{equation*}

\item \label{Minimizer sector}
Minimizer of Dirichlet integrals over sectors of annulus (ramp functions)\\
Let $Y:=\{g\in C^{1}(D_{\delta,1}(0,\theta_{0})): g(re^{0})=1$ and $g(re^{i\theta_{0}})=0$ $\forall r\in (\delta,1) \}$, for some constants $\delta<1,\theta_{0}\in [0,2\pi]$. Then there exists $\mu\in Y$ s.t. $\forall g\in Y$

\begin{equation*} I_{D_{\delta,1}(0,\theta_{0})}(g)\geq I_{D_{\delta,1}(0,\theta_{0})}(\mu)=\dfrac{log(\frac{1}{\delta})}{\theta_{0}} \end{equation*}
\end{enumerate}

\end{thm}

\subsection{Proof of Dubinin's Theorem}

\begin{proof}

For convenience let $\Omega:=\Omega_{\alpha}$, $S=S_{\alpha}$ and $\widetilde{\Omega}=\Omega_{\widetilde{\alpha}}$, $\widetilde{S}=S_{\widetilde{\alpha}}$. First, we simplify the problem. Consider origin-fixing conformal maps $F_{\Omega}:\Omega\to \mathbb{D}$ and $F_{\widetilde{\Omega}}:\widetilde{\Omega}\to \mathbb{D}$ (these maps can be obtained by extending the conformal map of a sector via the Schwartz reflection principle). We will show that $|F_{\Omega}(S)|\leq  |F_{\widetilde{\Omega}}(\widetilde{S})| \Leftrightarrow \omega(0,S,\Omega)\leq  \omega(0,\widetilde{S},\widetilde{\Omega})$. \\
Since $\omega(0,S,\Omega)=\omega(F_{\Omega}^{-1}(0),S,\Omega)$ is harmonic, the mean value property yields
\begin{align*}
\omega(0,S,\Omega)&=\frac{1}{2\pi}\int_{\partial \mathbb{D}}\omega(F_{\Omega}^{-1}(z),S,\Omega)dz\\
&=\frac{1}{2\pi}\int_{\partial \mathbb{D}}\chi_{F_{\Omega}(S)}dz\\
&=\frac{1}{2\pi}|F_{\Omega}(S)|.
\end{align*}
Thus, it suffices to show  $|F_{\Omega}(S)|\leq  |F_{\widetilde{\Omega}}(\widetilde{S})|$. Dubinin's idea was to make the slits coincide, so that the lengths of $S$ and $F_{\Omega}^{-1}(F_{\widetilde{\Omega}}(\widetilde{S})))$ can be compared (a priori F and $F_{\widetilde{\Omega}}$ act on slits at different positions). Here is the outline:\\

Step 1: Simplify the problem to the annulus $A(\delta,1)$ i.e. $|F_{ \widetilde{\Omega}(\delta)}(\widetilde{S})|\leq |F_{\Omega(\delta)}(S)|$.\\ 

Step 2: Consider minimizer $f$ of Dirichlet integral in $Y:=\{g\in C^{1}(D_{\delta,1}(0,\theta_{0})): g(re^{0})=1$ and $g(re^{i\theta_{0}})=0$ $\forall r\in (\delta,1) \}$ from \thmref{Dirichlet integral}.\ref{Minimizer sector}.\\

Step 3: Construct another minimizer $f_{4}\in Y$ in terms of conformal maps $F_{ \widetilde{\Omega}(\delta)},F_{\Omega(\delta)}(S)$ s.t. $I_{D_{\delta ,1}(0,\theta_{0})}(f_{4})\leq  I_{D_{\delta ,1}(0,\theta_{0})}(f)$ and conclude that $f_{4}$ cannot be identically 1 on $K=[a,1]$.\\

Step 4: Assume $|F_{\Omega}(S)|\geq  |F_{\widetilde{\Omega}}(\widetilde{S})|$, and obtain contradiction by showing that $f_{4}$ must be identically 1 on $K=[a,1]$.\\

Step 1:\\
Instead of $\Omega$, consider $\Omega(\delta):=\Omega\cap A(\delta,1)$ for $\delta\in (0,a)$. By the mapping theorem for doubly connected domains, there exists conformal maps $F_{\Omega(\delta)}:\Omega(\delta)\to A(\varepsilon,1)$, where $\varepsilon:=g(\delta)$ is a function of $\delta$ \cite{Duren_1983}. Then $Carath\acute{e}odory ~convergence~theorem$ yields $\lim\limits_{\delta \to 0}F_{\Omega(\delta)}= F_{\Omega}$ uniformly on each compact subset of $\Omega$. Thus, it suffices to prove
\begin{equation*}  |F_{ \widetilde{\Omega}(\delta)}(\widetilde{S})|\leq |F_{\Omega(\delta)}(S)|\end{equation*}

\begin{figure}[h]
\centering
\includegraphics[scale=0.3]{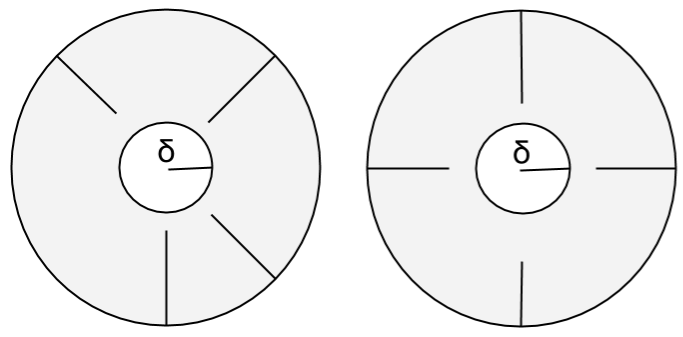}
\caption{Domains $\Omega(\delta)$ and $\widetilde{\Omega}(\delta)$ respectively}
\label{fig:my_label}
\end{figure}

We have the following lemma for annuli $A(\widetilde{\varepsilon},1)$, $A(\varepsilon,1)$, which we will need later.
\begin{lem}\label{epsilon}
For the above inner radii $\varepsilon$ and $\widetilde{\varepsilon}$, it holds that for all $\delta $
\begin{equation*} \varepsilon\leq \widetilde{\varepsilon}~\text{or equivalently }~A(\widetilde{\varepsilon},1)\subseteq A(\varepsilon,1) \end{equation*}
with strict inequality unless $\Omega$ can be obtained from $\widetilde{\Omega}$ by a rotation about the origin.
\end{lem}

\begin{proof}
First, we will get an inequality from \thmref{Dirichlet integral}.\ref{annulus integral} and then use the $2\pi (log(\frac{1}{\varepsilon}))^{-1}$ expression for harmonic functions on the annulus $A(\varepsilon,1)$. Let  $h:\Omega\to [0,1]$, $\widetilde{h}:\widetilde{\Omega}\to [0,1]$ be the functions defined by 
\begin{equation*} h(z):=\omega(z,|z|=\delta,\Omega(\delta))~\text{and}~\widetilde{h}(z):=\omega(z,|z|=\delta,\widetilde{\Omega}(\delta)) \end{equation*}

Let $\widetilde{h}^{dis}$ be a desymmetrization about $\alpha$. Since h is harmonic, we have by \thmref{Dirichlet integral}.\ref{harmonic dirichlet integral}
\begin{equation*}  I_{A(\delta,1)}(h) \leq   I_{A(\delta,1)}(\widetilde{h}^{dis}). \end{equation*}
Next we show that this implies $2\pi (log(\frac{1}{\varepsilon}))^{-1}\leq 2\pi (log(\frac{1}{\widetilde{\varepsilon}}))^{-1}$. Firstly, since desymmetrization preserves the Dirichlet integral ,
\begin{equation*}  I_{A(\delta,1)}(\widetilde{h})= I_{A(\delta,1)}(\widetilde{h}^{dis}). \end{equation*}
Secondly, $h(F_{\delta}^{-1}(z))$ is harmonic on the annulus $A(\varepsilon,1)$ with boundary values 0 on $\{|z|=1\}$ and 1 on $\{|z|=\varepsilon\}$. Similarly, for $\widetilde{h}(\widetilde{F}_{\delta}^{-1}(z))$ on $A(\widetilde{\varepsilon})$. Thus by conformal invariance of Dirichlet integral,
\begin{align*}
2\pi (log(\frac{1}{\varepsilon}))^{-1}&\stackrel{\thmref{Dirichlet integral}\ref{annulus integral}}{=}I_{A(\varepsilon,1)}(h(F_{\delta}^{-1}(z)))\\
&\stackrel{conformal}{=}I_{A(\delta,1)}(h)\\
&\leq I_{A(\delta,1)}(\widetilde{h}^{dis})=I_{A(\delta,1)}(\widetilde{h})\\
&=I_{A(\varepsilon,1)}(\widetilde{h}(\widetilde{F}_{\delta}^{-1}(z)))\\
&\stackrel{\thmref{Dirichlet integral}\ref{annulus integral}}{=}2\pi (log(\frac{1}{\widetilde{\varepsilon}}))^{-1}\\
&\Rightarrow \varepsilon\leq \widetilde{\varepsilon}.\\
\end{align*}
If $\Omega(\delta)$ and $\widetilde{\Omega}(\delta)$ are not rotation of each other, $\widetilde{h}^{dis}$ is not $C^{1}$ and thus not harmonic on $D_{\delta,1}(0,a_{j})$. Thus, the first inequality is strict  $I_{A(\delta,1)}(h) < I_{A(\delta,1)}(\widetilde{h}^{dis})$.

\end{proof}

Step 2:\\
Next we define on $A(\delta,1)$ a different auxiliary function f. From \thmref{Dirichlet integral}.\ref{Minimizer sector}, there exists minimizer $\mu$ on $D_{\delta,1}(0,\frac{\pi}{n})$. Then we extend this to function f on $A(\delta,1)$ by using the Schwartz reflection principle repeatedly i.e. we reflect $\mu$ across the lines $\widetilde{a}_{j}:=argz=\frac{j\pi}{n}$ for $j=1,...,2n-1$ and extend $\mu$ to function f defined on all of $A(\delta,1)$. 
Then for $r\in (\delta,1]$ it holds that $f(re^{i\widetilde{a}_{j}})=1\stackrel{\pref{desy prop}.\ref{equal on slits}}{\Rightarrow} S\subset\{z:f^{dis}(z)=1\}\Rightarrow f^{dis}(S)\subset\partial \mathbb{D}$.

Step 3:\\
The rest of the proof involves constructing a comparison function $f_{4}$ in terms of $F_{\Omega(\delta)}, F_{\widetilde{\Omega}(\delta)}$, with
\begin{equation*} I_{D_{\delta ,1}(0,\theta_{0})}(f_{4})\leq I_{D_{\delta ,1}(0,\frac{\pi}{n})}(f) 
\end{equation*}
This will be used to obtain a contradiction to $|F_{\Omega(\delta)}(S)|\geq |F_{ \widetilde{\Omega}(\delta)}(\widetilde{S})|$. Let $f_{1}:=f^{dis}$, $f_{2}:=f_{1}\circ F_{\Omega(\delta)}^{-1}$ and $f_{3}:A(\varepsilon,1)\to [0,1]$ be the n-fold circular symmetrization of $f_{2}$ i.e. $f_{3}(re^{i\theta}):=f_{2}^{*}(re^{in\theta})$\cite{Polya_Szego_1945}. Finally, let $f_{4}:=f_{3}\circ F_{\widetilde{\Omega}(\delta)}^{-1}$, which is well-defined since $A(\varepsilon,1)\subset A(\widetilde{\varepsilon},1)$ by \lemref{epsilon}.
\begin{lem}\label{f4_inequality}
For $f,f_{4}$ as above
\begin{equation*}  I_{D_{\delta ,1}(0,\frac{\pi}{n})}(f_{4})\leq I_{D_{\delta ,1}(0,\frac{\pi}{n})}(f)\end{equation*}
The inequality is strict unless the points $a_{j}$ are evenly spaced.
\end{lem}
\begin{proof}
The function f is symmetric on each sector since it was defined by reflecting it across $\widetilde{a}_{j}$; thus its Dirichlet integral over each sector is the same,
\begin{align*}
I_{D_{\delta ,1}(0,\frac{\pi}{n})}(f)&=\frac{1}{2n}I_{A(\delta,1)}(f)\\
&\stackrel{\pref{desy prop} }{=}\frac{1}{2n}I_{A(\delta,1)}(f_{1})\\
&\stackrel{conformal}{=}\frac{1}{2n}I_{A(\varepsilon,1)}(f_{2})\\
&\stackrel{\thmref{Dirichlet integral}.\ref{Modified Polya inequality}}{\geq}\frac{1}{2n}I_{A(\widetilde{\varepsilon},1)}(f_{3})
\end{align*}
The $f_{3}$ and $\widetilde{F}_{\delta}$ are symmetric about points $\widetilde{\alpha}_{j}$ on each circle $A(\varepsilon,1)$ and thus  $f_{4}=f_{3}\circ \widetilde{F}_{\delta}$ is also symmetric. Thus,\
\begin{align*}
I_{D_{\delta ,1}(0,\frac{\pi}{n})}(f_{4})&=\frac{1}{2n}I_{A(\delta,1)}(f_{4})\\
&=\frac{1}{2n}I_{A(\widetilde{\varepsilon},1)}(f_{3}).\\
&\leq I_{D_{\delta ,1}(0,\frac{\pi}{n})}(f).
\end{align*}
If $a_{j}$ are not evently spaced, then we have proper subset $A(\widetilde{\varepsilon},1)\subset A(\varepsilon,1)$. Also, the total variation of $f_{3}$ on each circle $\{|z|=r\}$ is 2n because n-fold circularly symmetric functions have n bumps (it goes around $[0,2\pi]$ n times). Thus,
\begin{equation*} \int_{A(\varepsilon,1)\setminus A(\widetilde{\varepsilon},1)} |\partial_{\theta}f_{3}|^{2}+|\partial_{r}f_{3}|^{2} r d\theta dr>0, \end{equation*}
and in turn $I_{A(\varepsilon,1)}(f_{2})>I_{A(\widetilde{\varepsilon},1)}(f_{3})$. The result follows.

\end{proof}

Step 4: \\
Finally we prove $|F_{\widetilde{\Omega}}(\widetilde{S})|\geq |F_{\Omega}(S)|$ via contradiction. Since $f=\mu$ has the smallest Dirichlet integral on sector $D_{\delta ,1}(0,\frac{\pi}{n})$ for functions in $\{g\in C^{1}(D_{\delta,1}(0,\frac{\pi}{n})): g(re^{0})=1$ and $g(re^{i\frac{\pi}{n}})=0$ $\forall r\in (\delta,1) \}$, the following lemma implies one of the following cases:
\begin{flalign*}
&f_{4}=f ~\text{on D,}\\
&f_{4} ~\text{is not admissible on D,}&\\
&f_{4} ~\text{is not identically 0 on}~ argz=\frac{\pi}{n},\text{ or,}&\\
&f_{4} ~\text{is not identically 1 on}~ argz=0.&
\end{flalign*}

Going through the definitions and using the equality case for \lemref{f4_inequality} shows that the first three cases are not true. Also, because $f((\delta,a))=1$, the maps' definitions yield that $f_{4}((\delta,a))=1$. Assuming, $(\delta,1)$ is on the real axis and using that $f_{4}$ is not identically 1 on $argz=0$, gives that $f_{4}$ cannot be identically 1 on $[a,1]$.\\
Assume $F_{\widetilde{\Omega}(\delta)}(\widetilde{S})\subset F_{\Omega(\delta)}(S)$, we will obtain a contradiction. First, we show that $f_{3}(z)=1$ for $z\in F_{\Omega(\delta)}(S)$. The measure $|F_{\Omega(\delta)}(S)|$ on $\partial \mathbb{D}$ is not affected by n-fold symmetrization because the $F_{\Omega(\delta)}(S)$ is split into n equal length intervals about the $\widetilde{a}_{j}$ on $\partial \mathbb{D}$. Thus, since $f_{1}=1$ on the slits of $\Omega(\delta)$, then $f_{3}(z)=f^{1}\circ F^{-1}_{\Omega(\delta)}(z)=1$ for $z\in F_{\Omega(\delta)}(S)$. \\
Since $[a,1]\subset \widetilde{S}\Rightarrow F_{\widetilde{\Omega}(\delta)}([a,1])\subset F_{\widetilde{\Omega}(\delta)}(\widetilde{S})\subset F_{\Omega(\delta)}(S)$, we get $f_{4}=f_{3}\circ (F_{\widetilde{\Omega}(\delta)})^{-1}$. This contradicts $f_{4}$ not being identically 1 on $[a,1]$. Hence $ F_{\Omega(\delta)}(S)\subset F_{\widetilde{\Omega}(\delta)}(\widetilde{S})$ on $\partial\mathbb{D}$.

\end{proof}

\subsection{Research Problems}
\begin{enumerate}
\item Let $g(\cdot, 0)$ denote the Green function of $\Omega_{\alpha}$ with pole at the origin and similarly $g^{*}(\cdot,0)$ Green function for $\Omega_{\widetilde{\alpha}}$. Does it hold that
\begin{equation*}  \int_{\Omega_{\widetilde{\alpha}} } g^{*}(y,0)dy\leq \int_{\Omega_{\alpha}}g(y,0)dy\end{equation*}
or in probabilistic interpretation
\begin{equation*} E_{0}T_{\Omega_{\widetilde{\alpha}} }\leq E_{0}T_{\Omega_{\alpha}}. \end{equation*}
\item Does the following stronger conjecture hold for $0<r<1$
\begin{equation*} \int_{|y|=r }g^{*}(re^{i\theta},0)d\theta\leq \int_{|y|=r }g(re^{i\theta},0)d\theta \end{equation*}
or in probabilistic interpretation for every $t>0$
\begin{equation*} P_{0}(T_{\Omega_{\widetilde{\alpha}}}>t)\leq P_{0}(T_{\Omega_{\alpha}}>t). \end{equation*}
\item For three slits (n=3 in Gonchar's problem) Dubinin's theorem \thmref{Gonchar-Dubinin theorem} follows from Baernstein's theorem: Let $\Omega_{\alpha}$ and $\Omega_{\widetilde{\alpha}}$ be as above and let $u(z):=\omega(z,S_{\alpha},\Omega_{\alpha})$ and $v(z):=\omega(0,S_{\widetilde{\alpha}},\Omega_{\widetilde{\alpha}})$ be their harmonic measures. Then for three slits (n=3), Baernstein showed\cite{Baernstein_1987}: Let $\Phi:\mathbb{R}\to \mathbb{R}$ be convex non-decreasing function $\Phi$ and $r\in (0,1)$ , then 
\begin{equation*} \int_{-\pi}^{\pi}\Phi(u(re^{i\theta}))d\theta \leq\int_{-\pi}^{\pi}\Phi(v(re^{i\theta}))d\theta . \end{equation*}
This inequality remains open for four or more slits (see \cite{Quine_1989,Betsakos_2001}).
\item For more general slit problems see \cite{Betsakos_2001} eg. disconnected intervals along the same slit.
\end{enumerate}

\pagebreak
\section{Extremal distance}
Extremal distance is a conformally invariant version of distance. As such it is a powerful tool for estimating conformal invariants like harmonic measure in terms of more geometric quantities. This approach leads to fundamental estimates for harmonic measure by extremal distance and the famous integral
\begin{eqnarray*}
\int_{\Omega} \frac{dx}{l(x)},
\end{eqnarray*}
where for domain $\Omega$ the $l(x)=|\Omega\cap \{z:Rez=x\}|$ is the height of a domain $\Omega$ at point x. In this section will exemplify this approach by using Carleman's method, . Baernstein says Carleman was the first to use integrals of the form $\int \frac{dx}{l(x)}$ to measure a domain and calls it "surely one of the most brilliant ideas in the history of complex function theory". 

\begin{thm}{Carleman} \cite{Marshall_Garnett_2005}\\
Let $\Omega\subset \mathbb{C}$ be a domain, $\Omega_{x}=\Omega\cap \{Rez=x\}$ and $E_{b}=\partial \Omega\cap \{Rez\geq b\}$. Suppose $\left |\Omega_{x}  \right |\leq M<\infty$ and let $l(x)$ denote the length of the longest interval in $\Omega_{x}$. Assume $z_{0}=x_{0}+iy_{0}\in B(z_{0},r_{0})\subset \Omega$. Then for $b>x_{0}$
\begin{equation*}   P_{z_{0}}(T_{E_{b}}=T_{\partial \Omega})\leq \dfrac{3M}{(2\pi r_{0}\int_{x_{0}}^{b}e^{2\pi \int_{x_{0}}^{t}\frac{dx}{l(x)}  }dt)^{\frac{1}{2}} }.\end{equation*}

\end{thm}

\begin{figure}[h]
\centering
\includegraphics[scale=0.4]{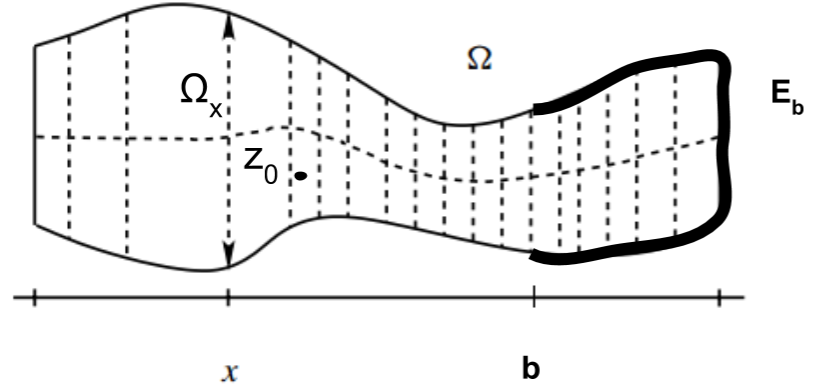}
\caption{Domain $\Omega$ and highlighted $E_{b}$ boundary}
\end{figure}

The proof is as follows: Carleman's idea was to find a differential inequality for the Dirichlet integral of the harmonic measure $\omega(z,E_{b},\Omega)$. Then by Green's theorem get an estimate for $\omega(z,E_{b},\Omega)$. To prove the theorem we may suppose that $\Omega$ is bounded, $\partial \Omega$ consists of finitely many analytic Jordan curves, and $inf\{Rez:z\in \Omega\}=0$. Write $\omega(z):=\omega(z,E_{b},\Omega)$ and define\\

\begin{equation*} A(t):=\int_{0}^{t}\int_{\Omega_{x}}|\triangledown \omega|^{2}dydx  . \end{equation*}

By the smoothness assumptions on $\partial \Omega$, the function $A(t)$ is continuously differentiable. Write $\Omega_{x}=\bigcup \Omega^{i}_{x}$ where $\{\Omega^{i}_{x}\}$ are the connected components of $\Omega_{x}$.\\

For the following lemma we will need Wirtinger's inequality: if g,g' are real valued continuous functions on the interval (a,b) and if $g(a)=g(b)=0$ then
\begin{equation*} \int_{a}^{b}(g')^{2}dx\geq (\frac{\pi}{b-a})^{2}\int_{a}^{b}g^{2}dx.  \end{equation*}
\begin{lem}{Carleman's differential inequality}\\
For $x\in (0,b)$, 
\begin{equation*} A'(x)\geq \frac{2\pi}{l(x)}A(x). \end{equation*}
\end{lem}
\begin{proof}
By the assumptions on $\partial \Omega$, for $x\in (0,b)$
\begin{equation*} A'(x)=\int_{\Omega_{x}}(\partial_{x}\omega)^{2}dy+\int_{\Omega_{x}}(\partial_{y}\omega)^{2}dy. \end{equation*}
Because $\omega=0$ on $\partial \Omega^{i}_{x}$, Wirtinger's inequality gives 
\begin{equation*}  \int_{\Omega^{i}_{x}}(\partial_{y}\omega)^{2}dy \geq (\frac{\pi}{|\Omega^{i}_{x}|})^{2} \int_{\Omega^{i}_{x}} \omega^{2}dy,\end{equation*}
and hence

\begin{equation*}\int_{\Omega_{x}}(\partial_{y}\omega)^{2}dy \geq (\frac{\pi}{|l(x)|})^{2} \int_{\Omega_{x}} \omega^{2}dy.  \end{equation*}

By Green's theorem
\begin{equation*} A(x)=\int_{\Omega_{x}}\omega \partial_{x}\omega dy, \end{equation*}

and by the Cauchy-Schwarz inequality

\begin{equation*} \int_{\Omega_{x}}(\partial_{x}\omega)^{2}dy \geq \frac{A^{2}(x)}{\int_{\Omega_{x}}\omega^{2}dy  }. \end{equation*}

Thus, in conclusion
\begin{equation*} A'\geq \frac{A^{2}(x)}{\int_{\Omega_{x}}\omega^{2}dy  }+(\frac{\pi}{|l(x)|})^{2} \int_{\Omega_{x}} \omega^{2}dy\geq \frac{2\pi}{l(x)}A(x). \end{equation*}

\end{proof}

\begin{proof}(Carleman's)
The Dirichlet integral A(t) is connected to the harmonic measure $\omega(z_{0})$ via the function $\phi(x)=\int_{\Omega_{x}}\omega^{2}dy$. By Harnack's inequality $\omega\geq \frac{\omega(z_{0})}{3}$ on $B(z_{0},\frac{r_{0}}{2})$, so that \\

\begin{equation*}\omega^{2}(z_{0})\leq \frac{9\phi(x_{0})}{r_{0}}.  \end{equation*}
On the other hand, because $\omega=0$ on $\partial \Omega$,
\begin{equation*} \phi'(x)=2\int_{\Omega_{x}}\omega \partial_{x}\omega dy. \end{equation*}
Because $A(x)=\int_{\Omega_{x}}\omega \partial_{x}\omega dy$ , $\phi'(x)=2A(x)$ and so Carleman's differential inequality reads
\begin{equation*} \frac{\phi''(x)}{\phi'(x)}\geq \frac{2\pi}{l(x)}. \end{equation*}
Now set $\mu(x)=\frac{2\pi}{l(x)}$ and $\psi(x)=\int_{0}^{x}e^{\int_{0}^{t}d\mu}dt$, so that 
\begin{equation*} \frac{\psi''}{\psi'}=\frac{2\pi}{l(x)}. \end{equation*}
Then the differential inequality can be rewritten as

\begin{equation*} (log\frac{\phi'}{\psi'})'=\frac{\phi''}{\phi'}-\frac{\psi''}{\psi'}\geq 0. \end{equation*}

Therefore $\frac{\phi'}{\psi'}$ is non-decreasing. Because $\psi'>0$, we obtain

\begin{equation*} \phi'(x)\psi'(t)\leq \phi'(t)\psi'(x) \end{equation*}
whenever $0<x<t$. Because $\phi(0)=\psi(0)=0$, integrating the above inequality from 0 to x gives 

\begin{equation*}  \phi(x)\psi'(t)\leq \phi'(t)\psi(x), \end{equation*}
and integrating again from x to t then gives\\

\begin{equation*} \phi(x)\psi(t)\leq \phi(t)\psi(x), \end{equation*}
whenever $0<x<t$. Increasing $\Omega$ in $\{Rez<x_{0}\}$ increases $\omega$ but does not change the right side of the inequality in the theorem, so we may assume that $\mu(x)=\frac{2\pi}{l(x)}=\frac{2\pi}{M}$ on $x<x_{0}$. Then $\psi(x_{0})=\frac{M}{2\pi}(e^{\frac{2\pi x_{0}}{M}}-1)$ and
\begin{equation*} \psi(b)=\psi(x_{0})+\int_{x_{0}}^{b}e^{\int_{0}^{x_{0}}\mu(s)ds  }e^{\int_{x_{0}}^{t}\mu(s)ds  }dt\geq \psi(x_{0})(1+\frac{2\pi}{M}\int_{x_{0}}^{b}e^{\int_{x_{0}}^{t}\mu(s)ds  }dt). \end{equation*}

Now $\phi(b)\leq |\Omega_{b}|\leq M$, so that by above 
\begin{equation*} \omega(z_{0})^{2}\leq \frac{9\phi(z_{0})}{r_{0}}\leq \frac{9M}{r_{0}}\frac{\psi(x_{0})}{\psi(b)}\leq \frac{9M}{r_{0}}(1+\frac{2\pi}{M}\int_{x_{0}}^{b}e^{\int_{x_{0}}^{t}\mu(s)ds  }dt)^{-1}. \end{equation*}

Thus, by taking square root
\begin{equation*}   P_{z_{0}}(T_{E_{b}}=T_{\partial \Omega})=\omega(z_{0})\leq \frac{3M}{(2\pi r_{0}\int_{x_{0}}^{b}e^{2\pi \int_{x_{0}}^{t}\frac{dx}{l(x)}  }dt)^{\frac{1}{2}} }.\end{equation*}

\end{proof}

\pagebreak

\section{Strong Markov property}\label{Strong Markov property}

The Markov property intuitively states that knowing the current position of a random process yields as much information as knowing the entire history of positions up to that point. In this section, we will exemplify uses of this property by proving the isoperimetric for exit probability under polarization. We will use these results in the next sections to obtain elementary proofs of the isoperimetrics for principal eigenvalues and capacities.

\subsection{Recursive Strong Markov property}

Let $D\subset \mathbb{R}^{n}$ be open with $n\geq 2$, $A\subset D^{c}$ be Borel set with $A^{\sigma}\subset (D^{\sigma})^{c}$. Then for $x\in D$ and $t>0$, the $p_{D}(t,x,A)$ denotes the probability that BM starting at x does not exit D for $s\leq t$ and $B_{t}\in A$.

\begin{thm}\label{exit probability polarization}
For notation as above it holds that for $x\in D$ and $t>0$
\begin{equation*} 
p_{D}(t,x,A)\leq p_{D^{\sigma}}(t,x^{\sigma},A^{\sigma}). 
\end{equation*}
\end{thm}

\begin{proof}

The proof is as follows: By repeatedly applying strong Markov property infinitely many times, we will express $p_{D}(t,x,A)$ only in terms of the exit probabilities from domains $D_{+},D_{-},D_{0}$, where the inequality is clear by symmetry and monotonicity domain for exit probability.\\
Let $D_{+}:=D\cap \mathbb{H}^{+}$, $D_{-}:=D\cap \mathbb{H}^{-}$ and $D_{0}:=D\cap \sigma D$.  We first assume that D is bounded and $\mathbb{R}\cap  \overline{(D_{+}\cap \partial  D_{0})}=\varnothing$ and $\mathbb{R}\cap \overline{(D_{-}\cap \partial  D_{0})}=\varnothing$; we will remove them after. For more details see \cite{Betsakos_1998}.\\
\begin{figure}
\includegraphics[scale=0.4]{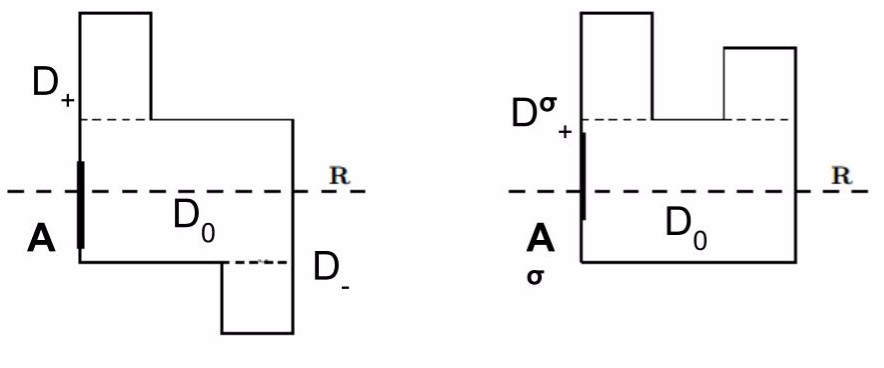}
\caption{Domain D and polarized domain $D^{\sigma}$}
\label{fig:my_label}
\end{figure}

Firstly, the hitting probability of A while avoiding $\partial D\setminus A$, splits into two terms: the h.p. of A while avoiding $\partial D_{0}\setminus A$ plus the h.p. of A while having hit $\partial D_{0}\setminus A$ but still avoiding $\partial D\setminus A$. The second term can be expressed using the Strong Markov property as\\
\begin{align*}
p_{D}(t,x,A)=&p_{D_{0}}(t,x,A)+\int_{0}^{\infty}\int_{D\cap \partial  D_{0}}P_{x}(B_{\tau_{D_{0}}}\in dl, \tau_{D_{0}}\in dt_{1}   )p_{D}(t-t_{1},l,A)\\
=&p_{D_{0}}(t,x,A)+\int_{0}^{\infty}\int_{D_{+}\cap \partial  D_{0}}P_{x}(B_{\tau_{D_{0}}}\in dl, \tau_{D_{0}}\in dt_{1}   )p_{D}(t-t_{1},l,A)\\
&+\int_{0}^{\infty}\int_{D_{-}\cap \partial  D_{0}}P_{x}(B_{\tau_{D_{0}}}\in ds, \tau_{D_{0}}\in dt_{1}   )p_{D}(t-t_{1},s,A)
\end{align*}

Secondly, applying it again for $p_{D}(t-t_{1},l,A)$ and $p_{D}(t-t_{1},s,A)$, yields
 \begin{equation*} p_{D}(t-t_{1},l,A)=p_{D}(t-t_{1},l,D_{+})+\int_{0}^{\infty}\int_{\mathbb{R}\cap D_{0}}P_{l}(B_{\tau_{D_{+}}}\in dr, \tau_{D_{+}}\in dt_{2}   )p_{D}(t-t_{1}-t_{2},r,A) \end{equation*}\\

\begin{equation*}p_{D}(t-t_{1},s,A)=p_{D}(t-t_{1},s,D_{+})+\int_{0}^{\infty}\int_{\mathbb{R}\cap D_{0}}P_{l}(B_{\tau_{D_{-}}}\in dr, \tau_{D_{-}}\in dt_{2}   )p_{D}(t-t_{1}-t_{2},r,A)   \end{equation*}

and so substituting them yields
\begin{align*}
&p_{D}(t,x,A)=p_{D_{0}}(t,x,A)+\int_{0}^{\infty}\int_{D_{+}\cap \partial  D_{0}}P_{x}(B_{\tau_{D_{0}}}\in dl, \tau_{D_{0}}\in dt_{1}   )p_{D_{+}}(t-t_{1},l,A)\\
&+\int_{0}^{\infty}\int_{0}^{\infty}\int_{D_{+}\cap \partial  D_{0}}\int_{D_{0}\cap \mathbb{R}}P_{x}(B_{\tau_{D_{0}}}\in dl, \tau_{D_{0}}\in dt_{1}   )P_{l}(B_{\tau_{D_{+}}}\in dr, \tau_{D_{+}}\in dt_{2}   )p_{D}(t-t_{1}-t_{2},r,A)\\
&+\int_{0}^{\infty}\int_{D_{-}\cap \partial  D_{0}}P_{x}(B_{\tau_{D_{0}}}\in dl, \tau_{D_{0}}\in dt_{1}   )p_{D_{-}}(t-t_{1},l,A)\\
&+\int_{0}^{\infty}\int_{0}^{\infty}\int_{D_{-}\cap \partial  D_{0}}\int_{D_{0}\cap \mathbb{R}}P_{x}(B_{\tau_{D_{0}}}\in ds, \tau_{D_{0}}\in dt_{1}   )P_{s}(B_{\tau_{D_{-}}}\in dr, \tau_{D_{+}}\in dt_{2}   )p_{D}(t-t_{1}-t_{2},r,A)
\end{align*}

We observe that terms of the form $p_{D}(\cdot,r,A)$ reappeared. Thus, we can repeat these two steps n-times to get an expression of $p_{D}(t,x,A)$ in terms of
\begin{equation*}
\begin{matrix}
 (1)p_{D_{+}}(\cdot, l,A,)  &  &(2)p_{D_{-}}(\cdot, l,A) &  & (3) p_{D}(\cdot,l,A )   \\ 
 (4) P_{l}(B_{\tau_{D_{+}}}\in ds, \tau_{D_{+}}\in dt   ) &  &  (5)P_{l}(B_{\tau_{D_{-}}}\in ds, \tau_{D_{-}}\in dt   )&  &   (6) P_{x}(B_{\tau_{D_{0}}}\in ds, \tau_{D_{0}}\in dt   )  
\end{matrix}
\end{equation*}
\\Similarly,  $p_{D^{\sigma}}(\cdot,r,A^{\sigma})$ has an expression in terms of\\
\begin{equation*}
\begin{matrix}
 (1)p_{(D^{\sigma})_{+}}(\cdot, l,A^{\sigma},)  &  &(2)p_{(D^{\sigma})_{-}}(\cdot, l,A^{\sigma}) &   (3)p_{D^{\sigma}}(\cdot,l,A^{\sigma} )     \\ 
 (4) P_{l}(B_{\tau_{(D^{\sigma})_{+}}}\in ds, \tau_{(D^{\sigma})_{+}}\in dt   )&  &  (5)P_{\sigma  l}(B_{\tau_{(D^{\sigma})_{+} }} \in ds, \tau_{(D^{\sigma})_{+}}\in dt   )&   (6) P_{x^{\sigma}}(B_{\tau_{(D^{\sigma})^{0}}}\in ds, \tau_{(D^{\sigma})^{0}}\in dt   )  
\end{matrix}
\end{equation*}
\\Because of domain monotonicity it holds that
\begin{align*}
 (1)p_{D_{+}}(\cdot, l,A,)&\leq p_{(D^{\sigma})_{+}}(\cdot, l,A^{\sigma},)\\
 (2)p_{D_{-}}(\cdot, l,A)&\leq p_{(D^{\sigma})_{-}}(\cdot, \sigma l,A^{\sigma}) \\
(4)P_{l}(B_{\tau_{D_{+}}}\in ds, \tau_{D_{+}}\in dt   )&\leq P_{l}(B_{\tau_{(D^{\sigma})_{+}}}\in ds, \tau_{(D^{\sigma})_{+}}\in dt   ) \\
(5)P_{l}(B_{\tau_{D_{-}}}\in ds, \tau_{D_{-}}\in dt   )&\leq P_{\sigma l}(B_{\tau_{(D^{\sigma})_{+}}}\in ds, \tau_{(D^{\sigma})_{+}}\in dt   ).
\end{align*}

Because of symmetry, $P_{x}(B_{\tau_{D_{0}}}\in ds, \tau_{D_{0}}\in dt   )=P_{x^{\sigma}}(B_{\tau_{D_{0}}}\in ds, \tau_{D_{0}}\in dt   )$. Finally, we will show that $p_{D}(\cdot,l,A )$ and $ p_{D^{\sigma}}(\cdot,l,A^{\sigma} )$ will vanish by repeating the argument above for $n\to \infty$. We will just do it for $p_{D}(\cdot,l,A )$.\\\\
At the nth-iteration, there will be $2^{n}$ integrals $\{I_{j}\}_{j=1}^{2^{n}}$ that contain $(3)$; these also contain n-factors of (4),(5),(6). Because of D being bounded and $\mathbb{R}\cap \overline{(D_{+}\cap \partial  D_{0})}=\varnothing$ and $\mathbb{R}\cap \overline{(D_{-}\cap \partial  D_{0})}=\varnothing$, there exists constant $M<1$ s.t.
 \begin{equation*}  P_{l}(B_{\tau_{D_{+}}}\in ds, \tau_{D_{+}}\in dt   )\leq M~\text{and}~P_{l}(B_{\tau_{D_{-}}}\in ds, \tau_{D_{-}}\in dt   )\leq M.\end{equation*}
In other words, the probability of immediately escaping from $D_{+},D_{-}$ is less than one. By symmetry 
\begin{equation*}  P_{x}(B_{\tau_{D_{0}}}\in ds, \tau_{D_{0}}\in dt   )\leq \frac{1}{2}.\end{equation*}
Therefore, $\sum_{j=1}^{2^{n}}I_{j}\leq \sum_{j=1}^{2^{n}}\frac{\delta^{n}}{2^{n}}=2\delta^{n}\to 0$ as $n\to \infty$. Similarly, integrals $I^{*}_{j}$ that contain $ p_{D^{\sigma}}(\cdot,l,A^{\sigma} )$ satisfy $\sum_{j=1}^{2^{n}}I^{*}_{j}\leq =2\delta^{n}\to 0$.\\

Next we remove the boundedness and assumption $\mathbb{R}\cap \overline{(D_{+}\cap \partial  D_{0})}=\varnothing$ and $\mathbb{R}\cap \overline{(D_{-}\cap \partial  D_{0})}=\varnothing$. The boundedness is removed by taking a sequence of increasing bounded open sets eg. $D_{m}:=D\cap \{|z|<m\}$ and using monotonicity of exit probabilities. Now onto removing the second assumption. For $n\in \mathbb{N}$, let
\begin{equation*} O_{n}:=\{z\in \mathbb{C}\setminus (A\cup \mathbb{R}):\inf\limits_{z\in \mathbb{R}\cap D_{0}}(z) <Re(z)<\sup \limits_{z\in \mathbb{R}\cap D_{0}}(z)~\text{and}~ -\frac{1}{n}<Im(z)<\frac{1}{n} \} \end{equation*}

and $D^{n}:=O_{n}\cup D$. Then $D^{n}$ satisfies $\mathbb{R}\cap \overline{(D^{n}_{+}\cap \partial  D^{n}_{0})}=\varnothing$ and $\mathbb{R}\cap \overline{(D^{n}_{-}\cap \partial  D^{n}_{0})}=\varnothing$. Thus, from above argument for $x\in \mathbb{R}\cap D^{n}_{0}$ 

\begin{equation*}  p_{D}(t,x,A)\leq p_{D^{\sigma}}(t,x^{\sigma},A^{\sigma}).\end{equation*}
Since $D_{n}$ is decreasing to D, the theorem follows.

\end{proof}

\subsection{Stopping times}
 The general strategy is to decompose the Brownian paths into appropriate stopping time segments.  For example, consider open $D\subset \mathbb{R}^{n}$ and subset $S\subset D$. Then we define inductively two types of sequences of stopping times for $k\geq 2$; the kth entry time and kth exit time respectively:
 \begin{equation*}\tau_{k,S}:=\inf\limits_{t>T_{k-1,S}}\{t\leq  T_{D}, X_{t}\in S\}~\text{ and}~ T_{k,S}:=\inf\limits_{t>\tau_{k,S}}\{t\leq  T_{D}, X_{t}\notin S\},  \end{equation*}

and for $k=1$\\
 \begin{equation*}  \tau_{1,S}:=\inf\limits_{t>0}\{t\leq  T_{D}, X_{t}\in S\}~\text{ and}~ T_{1,S}:=\inf\limits_{t>\tau_{1,S}}\{t\leq  T_{D}, X_{t}\notin S\}.\end{equation*}

Let $\mathbb{M}_{\kappa}^{n}$ denote for $\kappa=0$ the Euclidean space $\mathbb{R}^{n}$,for $\kappa=-1$ the hyperbolic $\mathbb{H}^{n}$ and for $\kappa=1$ the sphere $\mathbb{S}^{n}$. Setting $A=B$ in the theorem below, we get the difference term for the average of exit probabilities from $A$ and $A^{\sigma_{H}}$. This was proved in \cite{Burchard_Schmuckenschlager_2001}.\\

\begin{thm}\label{difference_term}
Let $A,B\subset \mathbb{M}^{n}_{\kappa}$ be Borel sets and $A^{\sigma_{H}},B^{\sigma_{H}}$ , then
\begin{equation*} \int_{B^{\sigma_{H}}}  P_{x}(T_{A^{\sigma_{H}}}>t) dx-\int_{B}  P_{x}(T_{A}>t) dx=\int_{B}\int_{A}  P_{x,y}(E^{t}) dxdy, \end{equation*}

where $E^{t}:=\{B_{[0,t]}(\omega): T_{A^{\sigma_{H}}}>t, m(B_{[0,t]}\cap (A\setminus \sigma_{H} A)\cap H^{+})>0~ \text{and}~ m(B_{[0,t]}\cap (\sigma_{H}A\setminus  A)\cap H^{+})>0\}$. In other words, it is the event that a Brownian path starting at x and ending at y does not leave $A^{\sigma}$ during $[0,t]$ and meets both $A\setminus \sigma_{H} A$ and $ \sigma_{H}A\setminus A$ during some subinterval where it remains in $H^{+}$.
\end{thm}
\begin{proof}
The proof is as follows: we decompose the Brownian paths into appropriate segments based on stopping times, where the result will follow by reflection.
\begin{align*}
\int_{B^{\sigma_{H}}}  P_{x}(T_{A^{\sigma_{H}}}>t) dx=&\int_{B}\int_{A}p_{A}(t,x,y)dydx\\
        =&\int_{H_{+}}\int_{H^{+}}p_{A}(t,x^{\sigma},y^{\sigma})1_{B}(x^{\sigma})1_{A}(y^{\sigma})\\
        +&p_{A}(t,x^{\sigma},\sigma y^{\sigma})1_{B}(x^{\sigma})1_{A}(\sigma y^{\sigma})\\
        +&p_{A}(t,\sigma x^{\sigma}, y^{\sigma})1_{B}(\sigma x^{\sigma})1_{A}(y^{\sigma})\\
        +&p_{A}(t,\sigma x^{\sigma},\sigma y^{\sigma})1_{B}(\sigma x^{\sigma})1_{A}(\sigma y^{\sigma}) dydx\\
        =:&\int_{H^{+}}\int_{H^{+}}\sum_{\pm}p_{A}(t,x^{\pm},y^{\pm})1_{B}(x^{\pm})1_{A}(y^{\pm}),\\
\text{where}~ x^{+}=x^{\sigma} \text{and } x^{-}=\sigma x^{\sigma}.
\end{align*}

So it suffices to show that
\begin{equation}\label{result}
\sum_{\pm}p_{A}(t,x^{\pm},y^{\pm})1_{B}(x^{\pm})1_{A}(y^{\pm})-\sum_{\pm}p_{A}(t,x^{\pm},y^{\pm})1_{B}(x^{\pm})1_{A}(y^{\pm})=\sum_{\pm}p_{E^{t}}(x^{\pm},y^{\pm})1_{B}(x^{\pm})1_{A}(y^{\pm})
\end{equation}

where $p_{E^{t}}(x^{\pm},y^{\pm})$ are paths in $E^{t}$ conditioned to start and end at $x^{\pm}$ and $y^{\pm}$ respectively. 
Next we construct an injective map L which assigns each Brownian path $B([0,t])$ contained in A to a corresponding path $LB([0,t])$ contained in $A^{\sigma}$. Let $K=A\cap\sigma A,$ $K^{+}= (A\sigma A)\cap H^{+}$ and $K^{-}= (A\sigma A)\cap H^{-}$. Also, define recursively a sequence of stopping times $\{T_{j}\}$ as follows:
First set, 
\begin{equation*}
T_{1}=\left\{\begin{matrix} 0 & ~\text{if}~ x\in B/ \sigma B\\ \inf\{0<s \leq min\{t,T_{A}\} \}& ,if~ x\in B\cap \sigma B\\ +\infty & else \end{matrix}\right.
\end{equation*}
Then given $T_{1},..., T_{j}$, set $T_{j+1}=t$ if
\begin{equation*}
\left\{\begin{matrix}B_{T_{j}}\in H^{+}, B_{s}\in A^{c}/ K^{-}~\text{for}~ T_{j}<s\leq t ~\text{and}~ B_{t}\in(A/ \sigma A)\cap H^{-}\\ B_{T_{j}}\in H^{-}, B_{s}\in A^{c}/ K^{+} ~\text{for}~ T_{j}<s\leq t ~\text{and}~ B_{t}\in(A/ \sigma A)\cap H^{+}\end{matrix}\right.
\end{equation*}

Otherwise set 
\begin{equation*}
T_{j+1}=\left\{\begin{matrix}inf\{T_{j}<s\leq t|B_{s}\in K^{-} \} & ~\text{if}~ B_{T_{j}}\in H^{+} \\ inf\{T_{j}<s\leq t|B_{s}\in K^{+} \} & ~\text{if}~ B_{T_{j}}\in H^{-} \end{matrix}\right.
\end{equation*}

\begin{figure}[h]
\centering
\includegraphics[scale=0.6]{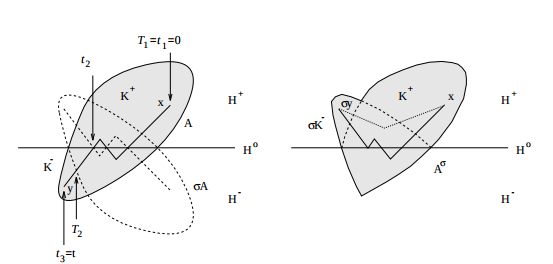}
\caption{Sets and Stopping times defined above. The dotted path on the right connecting x to $\sigma y$ is contained in $E^{t}$.}

\end{figure}

Denote by N the number of stopping times up to $min\{t,T_{A}\}$. Then, almost surely $N<\infty$ by the lemma below. 

\begin{lem}\label{finite}
Denote by N the number of stopping times up to $min\{t,T_{A}\}$. Then $N<\infty$. 
 
\end{lem}

\begin{proof}

Fix $n>1$, and choose an index set $J\subset \{1,...,n\}$. Let $S_{i}$ be the set of paths with $N=i$ as in the proof, and consider the image of $\bigcup_{i\geq n}S_{i}$ under the maps $L_{J}$ defined by 
\begin{equation*}
L_{J}(B_{s})=\left\{\begin{matrix}\sigma B_{s} ~\text{for}~ s\in [t_{j},t_{j+1}]&,if j\in J \\ B_{s} & else\end{matrix}\right.
\end{equation*}

Since the paths in $\bigcup_{i\geq n}S_{i}$ avoid $\sigma K\subset A^{\sigma}\setminus A$, the images of $L_{J}$ and $L_{J'}$ are disjoint for $J\neq J'$. In summary, for some constant $c>0$
 \begin{equation*} \sum_{\pm}P(\bigcup_{i\geq n}S_{i})\leq 2^{-n}\sum_{\pm}p(t,x^{\pm},y^{\pm})\leq \frac{c}{2^{n}}. \end{equation*}
\end{proof}

Let $S_{n}$ be the set of sample paths in A with $N=n$. By the intermediate value theorem and continuity of BM, every path $B_{s}$ in $S_{n}$ must hit $H^{0}$ at least once between $T_{j-1}$ and $T_{j}$. Let

 \begin{equation*}
 t_{j}=sup\{t\in (T_{j-1},T_{j})|B_{t}\in H^{0}\}~(j=2,...,n) 
 \end{equation*}

be the last time before $T_{j}$ that $B_{t}$ hits $H^{0}$, and set $t_{1}=0,t_{n+1}=min\{t,T_{A}\}$. Note that though $t_{j}$ is not a stopping time. For each path $B_{s}$ in $S_{n}$, the times $t_{j}$ cut the interval [0,t] into n subintervals, where $B_{s}$ hits $K^{+}$ and $K^{-}$ on alternating subintervals.\\
For $B_{s}\in S_{n}$, define 
\begin{equation*}
L(B_{s})=\left\{\begin{matrix}\sigma B_{s} ~\text{for}~ s\in [t_{j},t_{j+1}]&,\text{if}~ B_{T_{j}}\in K^{-} \\ B_{s} &\text{ else}\end{matrix}\right.
\end{equation*}

By construction, L maps a path in A which joins x with y on [0,t] to a path in $A^{\sigma}$ which joins x or $\sigma x$ with y or $\sigma y$. Since L is 1-1 and by reflection invariance of BM, the LHS of \eref{result} is nonnegative. By continuity, every path in the image of L meets $H^{0}$ on any interval where it hits both $K^{+}$ and $\sigma K^{-}$, so
 \begin{equation*} 
 LS_{n}\cap E^{t}=\varnothing 
 \end{equation*}
by definition of the event $E^{t}$. Conversely, if $E^{t}$ does not occur for a pth $B_{s}$ in $A^{\sigma}$, then one can construct an inverse image $L^{-1}B_{s}$ in A by reflecting the path on certain subintervals where it hits $\sigma K^{-}$. Then, the \eref{result} follows.
\end{proof}
\subsection{Research Problems}
For generalizations to Riemannian manifolds and open problems see \cite{Grigoryan_2002,Burchard_Schmuckenschlager_2001}.\\
(1) Is there a difference term for $ p_{D^{\sigma}}(t,x^{\sigma},A^{\sigma})-p_{D}(t,x,A)$?\\
(2) Is there an analogous result to \thmref{difference_term} for other symmetrizations eg. Steiner? Using repeated polarizations provides a difference term, but it would be interesting if there is an explicit probability event as with polarization.

\pagebreak

\section{Principal eigenvalue}\label{ch:Minimal_eigenvalue}

In this section, we will describe the connection between the principal eigenvalue of a domain $\Omega$ and the exit probability of BM from that domain. Further details can be found in \cite{Banuelos_2001,Schmuckenschlager_2011, Berg_Aqua_Sweers_2007}.
\subsection{Kac's formula}
Let $\Omega\subset \mathbb{R}^{n}$ be a bounded Lipschitz domain. By $\lambda_{k}, \phi_{k}$ we denote the kth eigenvalue and eigenfunction of the Dirichlet BVP, i.e.
\begin{equation*} -\frac{1}{2}\Delta \phi_{k}=\lambda_{k}\phi_{k}~\text{and}~\phi_{k}|_{\partial \Omega}=0. \end{equation*}
A probabilistic description of $\lambda_{1}(\Omega)$ is Kac's formula \cite{Port_Stone_1978} 
\begin{equation*} \lambda_{1}(\Omega)=-\lim\limits_{t\to \infty}\frac{1}{t}log P_{x}(T_{\Omega}>t). \end{equation*}

We prove this formula first. The density $p(t,x)=\frac{1}{(4\pi t)^{n/2}}e^{-\frac{|x|^{2}}{2t}}$ of BM killed at $\partial \Omega$ satisfies the heat equation over $\Omega$ with $p|_{\partial \Omega}=0$. Thus, it has an eigenfunction expansion i.e. $p(t,x)=\sum_{k=1}^{\infty}\left \langle p(0,x),\phi_{k}(x) \right \rangle_{H} e^{-\lambda_{k}t}\phi_{k}(x)$, eg. $H=L^{2}$. Thus, assuming $0<\lambda_{1}<\lambda_{2}\leq ...$, we obtain the formula from taking limit of
\begin{equation*}P_{p(0,x)}(T_{\Omega}>t)=\int_{\Omega}p(t,x)dx=e^{-\lambda_{1}t}[\left \langle p(0,x),\phi_{1}(x) \right \rangle_{H}\int \phi_{1}(x)dx+  O(e^{-(\lambda_{2}+\lambda_{1})t})].  \end{equation*}

\subsection{Faber-Krahn inequality and other results }

Thus, because of $P_{x}(T_{\Omega}>t)\leq P_{0}(T_{\Omega^{*}}>t)$ (proved in the "Brascramp-Lieb-Luttinger" section), we obtain 

\begin{thm}[Faber-Krahn inequality]
Let $\Omega\subset \mathbb{R}^{n}$ be a bounded Lipschitz domain then
\begin{equation*}\lambda_{1}(\Omega)\geq \lambda_{1}(\Omega^{*}). \end{equation*}
\end{thm}                                                 
It also follows that if $\lambda_{1}>0$, then $E_{x}(T_{\Omega})<\infty$: By the tail formulation with constant c 
\begin{equation*} E_{x}(T_{\Omega})=\int_{0}^{\infty}P_{x}(T_{\Omega}>t)dt=\int_{0}^{\infty}e^{-\lambda_{1}t}c +O(e^{-(\lambda_{1}+\lambda_{2})t})dt<\infty.. \end{equation*}

We also obtain a Brunn-Minkowski type inequality for principal eigenvalues.
\begin{thm}
Let $B,D\subset \mathbb{R}^{n}$ be balls , then 
\begin{equation*}\lambda_{1}(\frac{1}{2}(B+D))\leq \frac{1}{2}(\lambda_{1}(B)+\lambda_{1}(D))~\text{and}~\lambda_{1}(B\cap D)\leq \lambda_{1}(B)+\lambda_{1}(D).  \end{equation*}

\end{thm}

\begin{proof}

Before we proceed, we need the following interpolation inequality, which is interesting in itself.
\begin{lem}[Interpolation inequality]
Let $A,B,C\subset \mathbb{R}^{n} $ be open and $(1-\lambda)A+\lambda B\subset C$ for some $\lambda\in (0,1)$, then
\begin{equation*}P_{(1-\lambda)x+\lambda y}(T_{C}>t)\geq P_{x}(T_{A}>t)^{1-\lambda} P_{y}(T_{B}>t)^{\lambda}.  \end{equation*}
\end{lem}
\begin{proof}

Let $0<s_{1}<...<s_{m}\leq t$ be a finite subset of $(0,t]$. Notice that we can rewrite
\begin{equation*} P_{x}(\forall 1\leq j\leq m; B_{s_{j} }\in A)=P_{0}(Y\in A^{m}-x^{m}), \end{equation*}

where $Y=(X_{s_{1}},...,X_{s_{m}})$, $A^{m}=A\times \cdots\times A$ and $x^{m}=(x,...,x)$; they are all subset of $\mathbb{R}^{nm}$. Finally,
since $(1-\lambda)(A^{m}-x^{m})+\lambda(B^{m}-y^{m})\subset C^{m}-(1-\lambda)x^{m}-\lambda y^{m}$, the result follows from the log-concavity of the image Gaussian measure $\mu(A)=P_{0}(Y\in A^{m}-x^{m})$.
\end{proof}

The inequalities will follow from Kac's formula and the above lemma. For the first inequality,
\begin{align*}
\frac{1}{2}(\lambda_{1}(B)+\lambda_{1}(D))&=-[\lim\limits_{t\to \infty}\frac{1}{t}log P_{x}(T_{B}>t)^{\frac{1}{2}}+\lim\limits_{t\to \infty}\frac{1}{t}log P_{x}(T_{D}>t)^{\frac{1}{2}}]\\
&=-[\lim\limits_{t\to \infty}\frac{1}{t}log (P_{x}(T_{B}>t)^{\frac{1}{2}}P_{x}(T_{D}>t)^{\frac{1}{2}})]\\
&\geq -[\lim\limits_{t\to \infty}\frac{1}{t}log P_{\frac{1}{2}(x+y)}(T_{B+D}>t)]\\
&=\geq -[\lim\limits_{t\to \infty}\frac{1}{t}log P_{x+y}(T_{\frac{1}{2}(B+D)}>t)]\\
&=\lambda_{1}(\frac{1}{2}(B+D)).
\end{align*}
For the second inequality, by Lieb's result \cite{Lieb_1983}, there exists $x\in \mathbb{R}^{n}$ s.t.
\begin{equation*} \lambda_{1}(B\cap (D+x))\leq \lambda_{1}(B)+\lambda_{1}(D). \end{equation*}

Thus, by Kac's formula, it suffices to prove that $P_{0}(T_{B\cap (D+x)>t})\leq P_{0}(T_{B\cap D}>t)$. For any $x\in \mathbb{R}^{n}$ we have 
\begin{equation*} \frac{1}{2}(B\cap (D+x)+B\cap (D-x))\subset B\cap D~\text{and by symmetry }~P_{0}(T_{B\cap (D-x)}>t)=P_{0}(T_{B\cap (D+x)}>t). \end{equation*}

Thus, the result follows by above lemma 
\begin{equation*} P_{0}(T_{B\cap D}>t)\geq P_{0}(T_{B\cap (D+x)}>t)^{1/2} P_{0}(T_{B\cap (D-x)}>t)^{1/2}=P_{0}(T_{B\cap (D+x)}>t). \end{equation*}
\end{proof}

\subsection{Research Problems}
The following problems are from \cite{Henrot_2004}.\\
(1) Prove that the regular n-gone has the least first eigenvalue among all the n-gone of given area for $n\geq 5$.\\
(2) Let $\Omega$ be fixed domain and $B_{0}$ ball of fixed radius. Prove that $\lambda_{1}(\Omega\setminus B_{0})$ is minimal when $B_{0}$ intersects $\partial \Omega$ at a point and is maximum when $B_{0}$ is centered at a particular point of $\Omega$.\\
(3) Prove that a convex domain $\Omega^{*}\subset \mathbb{C}$ which minimizes $\lambda_{2}$ (among convex domains of given area) has two perpendicular axes of symmetry.\\
(4) In dimensions 2 and 3 , prove that the optimal domain for $\lambda_{3}$ is a ball and in dimension $n\geq 4$ the union of three identical balls.

\pagebreak

\section{Unified Symmetrization by Polarization}

Another approach to symmetrization can be based on polarization. This simple rearrangement was introduced by Wolontis (planar sets- 1952)\cite{Wolontis_1952}, Baernstein and Taylor (functions-1976) \cite{Wolontis_1952,Baernstein_Taylor_1976}. In \cite{Baernstein_Taylor_1976} Baernstein and Taylor provide easier proofs of some integral inequalities by replacing symmetrization with polarization. Furthermore, Dubinin and Solynin used it to prove inequalities for capacities. For more details and results see \cite{Solynin_Brock_2000, Schaftingen_2006}. In this section we will exemplify this method by proving the isoperimetric for capacities\cite{Betsakos_2004}.

\subsection{Approaching symmetrizations by polarizations}

Let $K\subset \mathbb{R}^{n}$ be a compact set with finite volume. Consider Steiner symmetrization of K wrt to hyperplane $H$, then $St(K)=St(K^{\sigma_{H}})=(St(K))^{\sigma_{H}}$. Thus, $d_{Ha}(K^{\sigma_{H}},St(K))< d_{Ha}(K,St(K))$, where $d_{Ha}$ is the Hausdorff distance. For repeated Steiner symmetrizations one can show \cite{Solynin_Brock_2000}:

\begin{prop}\label{polarization to steiner}
Let $K\subset \mathbb{R}^{n}$ be a compact set with finite volume, then there exists a sequence of polarizations $\{\sigma_{k}\}_{k\in \mathbb{N}}$ s.t.
\begin{equation*}\lim\limits_{k\to \infty}d_{Ha}(K^{\sigma_{1}\cdots \sigma_{k}}, St(K) )=0.  \end{equation*}

\end{prop}

This has the following corrolary. From the above convergence, given $r>0$, there exists $N\in \mathbb{N}$ s.t. $\forall k\geq N$
\begin{equation*}St(K)\subset K^{\sigma_{1}\cdots\sigma_{k}} +r\mathbb{B}^{n},  \end{equation*}

where $\mathbb{B}^{n}$ is the unit ball. Thus, since polarization is a smoothing transformation i.e. $A^{\sigma}+r\mathbb{B}^{n}\subset (A+r\mathbb{B}^{n})^{\sigma}$,
\begin{equation*} St(K)\subset K^{\sigma_{1}\cdots\sigma_{k}} +r\mathbb{B}^{n}\subset (K+r\mathbb{B}^{n})^{\sigma_{1}\cdots\sigma_{k}}. \end{equation*}

This decreasing chain along with the monotonicity of capacity will be used to obtain the isoperimetric result below.

\subsection{$\alpha$-Capacity decreases under Steiner symmetrization}
We already proved this result in section 1. We briefly repeat the definitions and theorem. The $\alpha$-Riesz kernel is
\begin{equation*} k_{\alpha}(x-y)=\frac{\Gamma(n-\frac{\alpha}{2})}{\Gamma(\frac{\alpha}{2}) \pi^{\frac{n}{2} 2^{\alpha-1}  }}\frac{1}{|x-y|^{n-\alpha}} , \end{equation*}
where $n\geq 2$ and $0<\alpha<n$. Let K be a compact set in $\mathbb{R}^{n}$, the $\alpha-$Riesz capacity of K is defined by\\
\begin{equation*} Cap_{\alpha}(K):=[\inf\limits_{\mu}\int\int k_{\alpha}(x-y) d\mu(x)d\mu(y)  ]^{-1}, \end{equation*}
where the infimum is taken over all probability Borel measures supported in K. If $\alpha=2$ and for $n\geq 3$,  this is the Newtonian capacity. 
\begin{thm}
Let $K\subset \mathbb{R}^{n}$ be a compact set with finite positive volume, then for $\alpha\in (0,2)$ and Steiner symmetrization
\begin{equation*} Cap_{\alpha}(K)\geq Cap_{\alpha}(St(K))\geq Cap_{\alpha}(K^{*}). \end{equation*}
\end{thm}
We start with the following lemma, which is of independent interest.
\begin{lem}\label{lemma symmetrization}
Let $K\subset \mathbb{R}^{n}$ be a compact set with finite positive volume, then for $\alpha\in (0,2)$
\begin{equation*} Cap_{\alpha}(K)\geq Cap_{\alpha}K^{\sigma_{H}} \end{equation*}
\end{lem}

\begin{proof}
Let $G_{\alpha,K^{c}}(x,y)$ be the Green function of $K^{c}$, Port proved the following identities between \cite{Port_Stone_1978}

\begin{equation*} Cap_{\alpha}(K)=lim_{|x|\to \infty}\frac{1}{c(n,\alpha)}|x|^{n-\alpha}P_{x}(T_{K^{c}}<\infty) \end{equation*}

and for $x\in K^{c}$
\begin{equation*}  1-P_{x}(T_{K^{c}}<\infty)=lim_{|y|\to \infty}\frac{1}{c(n,\alpha)}|y|^{n-\alpha}G_{\alpha,K^{c}}(x,y),\end{equation*}

where $c(n,\alpha):=\frac{\Gamma(\frac{n-\alpha}{2})}{\Gamma(|\frac{\alpha}{2}|)2^{\alpha}\pi^{\frac{n}{2}}}$. We have $G_{\alpha,K^{c}}(x,y)=\int_{0}^{\infty}p_{K^{c}}(t,x,y)dt$, where $p_{K^{c}}(t,x,y)$ is the transition probability from x to y at time t of the process $X_{t}$ killed on exiting $K^{c}$. \\

In "Strong Markov property" section we show that $p_{K^{c}}(t,x,y)\leq p_{(K^{c})^{\sigma_{H}}}(t,x^{\sigma_{H}},y^{\sigma_{H}})$ and so for $G_{\alpha,K^{c}}(x,y)=\int_{0}^{\infty}p_{K^{c}}(t,x,y)dt$ with $x,y\in K^{c}$, it holds that
\begin{equation*}G_{\alpha,K^{c}}(x,y)\leq G_{\alpha,(K^{\sigma_{H}})^{c} }(x^{\sigma_{H}},y^{\sigma_{H}}).  \end{equation*}
Thus, by the second identity above for $x\in H\cap K^{c}$ it holds that
\begin{equation*}P_{x}(T_{(K^{\sigma_{H}})^{c}  }<\infty)\leq P_{x}(T_{K^{c}}<\infty) ,  \end{equation*}

and so by the first identity
\begin{equation*} Cap_{\alpha}(K)\geq Cap_{\alpha}K^{\sigma_{H}}. \end{equation*}

\end{proof}

\begin{proof}

Let $K_{j}:=K+\frac{1}{j}\mathbb{B}^{n}$ and $\varepsilon>0$. The sequence $\{K_{j}\}$ is decreasing, and $\bigcap_{k=1}^{\infty}K_{j}=K$. Thus,  the continuity from above of Newtonian capacity (for this and other properties \cite{ Chung_1982, Landkoff_1972}) yields
\begin{equation*}  \lim\limits_{j\to \infty}Cap_{\alpha}K_{j}=Cap_{\alpha}K\end{equation*}
or equivalently there exists a $N_{1}\in \mathbb{N}$ s.t. 
\begin{equation*} Cap_{\alpha}(K)+\varepsilon\geq Cap_{\alpha}(K+\frac{1}{N_{1}}\mathbb{B}^{n}). \end{equation*}
Therefore, by discussion under \pref{polarization to steiner}, we can choose  $N_{2}\in \mathbb{N}$ s.t.
\begin{equation*} St(K)\subset (K+\frac{1}{N_{1}}\mathbb{B}^{n})^{\sigma_{1}\cdots\sigma_{N_{2}}}  \end{equation*}
and the monotonicity of $\alpha$-Riesz capacity yields
\begin{equation*} Cap_{\alpha}(St(K))\leq Cap_{\alpha}(K+\frac{1}{N_{1}}\mathbb{B}^{n})^{\sigma_{1}\cdots\sigma_{N_{2}}}. \end{equation*}

Finally, by \lemref{lemma symmetrization}
\begin{equation*} Cap_{\alpha}(K)+\varepsilon\geq Cap_{\alpha}(K+\frac{1}{N_{1}}\mathbb{B}^{n})\geq Cap_{\alpha}(K+\frac{1}{N_{1}}\mathbb{B}^{n})^{\sigma_{1}\cdots\sigma_{N_{2}}}\geq Cap_{\alpha}(St(K)) . \end{equation*}
By repeating this argument for a sequence of Steiner symmetrizations of K converging to $K^{*}$, gives the second inequality.
\end{proof}

\subsection{Research Problems}
\begin{enumerate}
\item Let $K\subset \mathbb{R}^{n}$ for $n\geq 3$ be a compact set with finite positive volume, then for $\alpha\in (2,n)$
\begin{equation*} Cap_{\alpha}(K)\geq Cap_{\alpha}. \end{equation*}
\item P$\acute{o}$lya-Szeg$\ddot{o}$ conjecture for electrostatic capacities: Let $K:=\{\Omega\subset \mathbb{R}^{3}:\Omega$  bounded and convex, $\mathcal{H}^{2}(\Omega)>0\}$, where $\mathcal{H}^{2}$ is the Hausdorff measure, then \cite{Fragala_2011} the conjecture is that the minimizer in K is attained by a disk
\begin{equation*} \sqrt{\frac{4\pi}{Area(\Omega)}}Cap(\Omega)\geq \inf\limits_{\Omega\in K} \sqrt{\frac{4\pi}{Area(\Omega)}}Cap(\Omega)=\sqrt{\frac{4\pi}{Area(\mathbb{D})}}Cap(\mathbb{D})= \frac{2\sqrt{2}}{\pi}\approx 0.9. \end{equation*}
\end{enumerate}

\pagebreak

\section{From the Sphere to the Euclidean space}

This section is about obtaining inequalities on the sphere and then transfering them to the Euclidean space by projecting. To exemplify this we will prove the isoperimetric for the Wiener sausage \cite{Ito_McKean_1965}. For any open set $A\subset \mathbb{R}^{d}$ with finite volume, Kesten  proved for $d\geq 3$ that\cite{Ito_McKean_1965}
\begin{equation*}\lim\limits_{t\to \infty}\frac{1}{t}E(vol(\bigcup_{s\leq t} B_{s}+A  ))=Cap(A),  \end{equation*}
where $Cap(A)$ is the Newtonian capacity of A. The quantity $E(vol(\bigcup_{s\leq t} B_{s}+A  )$ is called the \textit{Wiener sausage}. Various asymptotic results have been obtained about it \cite{Spitzer_1964, Donsker_Varadhan_1975, Gall_1988}.
\begin{figure}[h]
\centering
\includegraphics[scale=0.5]{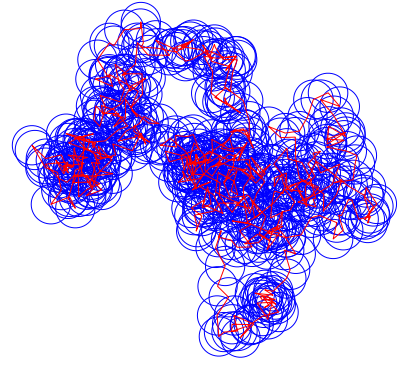}
\caption{Wiener sausage where A is a disk. The red line corresponds to the BM.}
\end{figure}

We will prove the following more general isoperimetric theorem , where we can allow A to vary with time. The isoperimetric result for capacities follows by setting $A_{s}=A$ for all $s\geq 0$.\\ 

\begin{thm}[Wiener sausage isoperimetric theorem]
Let $(A_{s})_{s\in \mathbb{R}^{+}}$ be open sets $\mathbb{R}^{d}$ and balls $(B_{x}(r_{s}))_{s\in \mathbb{R}^{+}}$ s.t. $vol(A_{s})=vol(B_{x}(r_{s}))$, then \cite{Peres_Sousi_2012}

\begin{equation*} E(vol(\bigcup_{s\leq t} B_{s}+A_{s}  )\geq E(vol(\bigcup_{s\leq t} B_{B_{s}}(r_{s}) ). \end{equation*}

\end{thm}

\subsection{Proof of the Wiener sausage isoperimetric theorem}
We only prove the analogous result below for finitely many open sets and random walk, in order to convey the main idea. Then, the general theorem follows from considering countable dense sets eg. dyadic cubes and Donsker's invariance principle (for more details see \cite{Peres_Sousi_2012}). 

\begin{prop}[Finitely many open sets]\label{finite case}
Let $\varepsilon>0$ and $(z_{*}(k))_{k\geq 0} $ be a random walk in $\mathbb{R}^{d}$ with $z_{*}(0)=0$ and transition kernel
\begin{equation*}p(x,y)=\frac{1}{\varepsilon^{n}\omega_{n}}1(\left \| x-y \right \|<\varepsilon).  \end{equation*}
Then. for any collection of open sets $(U_{k})_{k\geq 0}$ in $\mathbb{R}^{d}$, it holds for any $n\in \mathbb{N}$
\begin{equation*} E(vol(\bigcup_{k=0}^{n} z_{*}(k)+ )\geq E(vol(\bigcup_{k=0}^{n} B_{z_{*}(k)}(r_{k}) ). \end{equation*}
\end{prop}

\begin{proof}
The proof is as follows: starting from an isoperimetric result for random walks starting uniformlyon the sphere, we make precise the projecting of this inequality to the Euclidean space and using the mean of uniform distribution obtain the Wiener sausage inequality. We denote by $\mathbb{S}^{d}_{R}$ the sphere of radius R centered at 0, $\mu$ the surface area measure on the spere , $\rho(x,y)$ the geodesic distance between x,y and the geodesic cap centered at x of geodesic radius $\varepsilon$
\begin{equation*} C_{\varepsilon}(x):=\{z\in \mathbb{S}^{d}_{R}:\rho(x,z)<\varepsilon \}  . \end{equation*}
Let $\varepsilon>0$ and $\widetilde{\zeta}$ be a random walk on the sphere that starts from uniform point on the surface of the sphere $\mathbb{S}^{d}_{R}$ and has transition kernel given by 
\begin{equation*}\psi(x,y):=\frac{1}{\mu(C_{\varepsilon} (x))}1(\rho(x,y)<\varepsilon).  \end{equation*}
We call such a random walk the \textit{$\varepsilon-$cap walk}. For collection $(\Theta_{k})_{k\geq 0}$ Borel subsets of $\mathbb{S}^{d}_{R}$ and $\Theta:=\bigcup_{k=0}^{\infty}\Theta_{k}$ we define stopping time
\begin{equation*}\tau^{\Theta}:=min\{k\geq 0: \widetilde{\zeta}(k)\in \Theta_{k}\}     \end{equation*}

Then the following result follows from the BLL on the sphere \cite{Burchard_Schmuckenschlager_2001} (stated in the Brascamp-Lieb-Luttinger section)
\begin{lem}[Hitting probability isoperimetric on sphere]\label{isoperimetric}

Let $(\Theta_{k})_{k\geq 0}$ be Borel subsets of $\mathbb{S}^{d}_{R}$, then for all n it holds that
\begin{equation*}P(\tau_{\Theta}>n)\leq P(\tau_{C}>n)  \end{equation*}

where $C:=\bigcup_{k=0}^{\infty}C_{k}$ and $C_{k}$ is the geodesic cap centered at $(0,....,0,-R)$ such that its spherical measure $\mu(C_{k})=\mu(\Theta_{k})$.
\end{lem}

It is intuitive that as $\lim\limits_{R\to \infty}\mathbb{S}_{R}^{d}=\mathbb{R}^{d}$. The following lemmas make this precise in the situation we need for the proof of \thmref{finite case}. 

\begin{lem}\label{sphere}[From Sphere to Euclidean space]\\
Let $\pi:\mathbb{S}_{R}^{d}\to \mathbb{R}^{d}$ be the projection map $\pi(x_{1},...,x_{d+1})=(x_{1},...,x_{d})$ with inverse $\pi^{-1}((x_{1},...,x_{d}))=(x_{1},...,x_{d+1})\in \mathbb{S}_{R}^{d}$ s.t. $x_{d+1}\leq 0$. 
\begin{enumerate}

\item \label{Intermediate value property}[Intermediate value property]
Let $r,K>0$; then given $\delta>0$ there exists $R_{0}$ s.t. for all $x\in B_{K}(0)\subset\mathbb{R}^{n}$ it holds that

\begin{equation*} B_{r-\delta}(x)\subset \pi(C_{r}(\pi^{-1}(x)))\subset B_{r+\delta}(x)~\text{for all }~R>R_{0} \end{equation*}

\item \label{preserving measure} [Preserving measure] 
Let $K>0$ and $ A\subset \pi^{-1}(B_{K}(0))$ then
\begin{equation*} lim_{R\to \infty}\mu(A)=vol(\pi(A)). \end{equation*}

To \textit{couple} two rvs X ,Y from $(\Omega_{1},F_{1},P_{1})$,$(\Omega_{22},F_{2},P_{2})$ respectively; means to create a new probability space $(\Omega,F,P)$ over which there are rvs $X', Y'$  s.t. $X'\stackrel{d}{=}X$ and $Y'\stackrel{d}{=}Y$. The simplest way to couple two random walks is to force them to walk together. On every step, if X walks, so does Y. Thus, the difference between the two particles stays fixed. The following makes this precise in our context

\item\label{coupling} [Coupling random walks on sphere and Euclidean space]
Let $n,L>0$, $z$ an $\varepsilon-$ball walk in $\mathbb{R}^{d}$ started from a uniform point in $B_{L+n\varepsilon}(0)$, $\zeta$ an $\varepsilon-$ball walk in $\mathbb{S}^{d}_{R}$ started from a uniform point in $C(L):=\pi^{-1}(B_{L+n\varepsilon}(0))$. Then there exists a coupling of $z$ and $\pi(\zeta)$ s.t.
\begin{equation*} lim_{R\to \infty}P(\pi(\zeta_{k} )=z_{k}, \forall 0\leq k\leq m)=1.\end{equation*}

\end{enumerate}
\end{lem}

We start by putting the problem in the lemmas' context. We will assume $(U_{k})_{k\geq 0}$ are bounded since then we obtain the result by truncation. Thus, there exists large enough L s.t. $\bigcup_{i=0}^{n}U_{i}\subset B_{L}(0)$.  Next consider large enough radius R of $B_{0}(R)$ so that for $\varepsilon>0$

\begin{equation*} C(L)=\pi^{-1}(B_{L+n\varepsilon}(0)), \end{equation*}
where $C(L)$ is a cap centered at $(0,...,0,-R)$ of geodesic radius bigger than $L+n\varepsilon$. Finally, let $C_{k}$ be a geodesic cap centered at $(0,...,0,-R)$ s.t. $\mu(C_{k})=\mu(\pi^{-1}(U_{k}))$. \\
First, we prove 
\begin{equation*} P(\forall k=0,...,n, \pi(\zeta(k))\notin U_{k}  )\leq P(\forall k=0,...,n, \pi(\zeta(k))\notin \pi(C_{k})  ). \end{equation*}
From \lemref{isoperimetric} it holds that $P(\tau^{\pi^{-1}(U) }>n)\leq P(\tau^{C }>n)$. Also, by the Strong Markov property $P(\tau^{\pi^{-1}(U) }\leq n)=P(\widetilde{\zeta}(0)\in C(L) )P_{C(L)}(\tau^{\pi^{-1}(U) }\leq n)$, where $P_{C(L)}$ means starting uniformly from $C(L)$. Similarly, $P(\tau^{C} \leq n)=P(\widetilde{\zeta}(0)\in C(L) )P_{C(L)}(\tau^{C }\leq n)$. Thus, the inequality follows.

Secondly, let $\sigma_{U}:=min\{k\geq 0: z(k)\in U_{k}\}$ when $U:=\bigcup_{k} U_{k}$, then we prove that 
\begin{equation*}  
P(\sigma_{\pi(C)}>n)\geq P(\forall k, z(k)\notin U_{k})-2P(\exists k ~\text{ s.t. }~\pi(\zeta(k))\neq z(k) ).
\end{equation*}

We will denote $P(\text{coupling fails}):=P(\exists k$ s.t. $\pi(\zeta(k))\neq z(k) )$. By de Morgan's law
\begin{align*}
P(\bigcap_{k}^{n}z_{k}\notin U_{k})-P(\text{coupling fails})&=1-P(\bigcup_{k}^{n}z_{k}\in U_{k})-P(\text{coupling fails})\\
&\leq 1-P(\bigcup_{k}^{n}z_{k}\in U_{k}\cup \text{coupling fails})\\
&\leq P(\bigcap_{k}^{n}\pi(\zeta_{k})\notin U_{k}).
\end{align*}
Similarly, $P(\bigcap_{k}^{n}\pi(\zeta_{k})\notin \pi(C_{k}))-P(\text{coupling fails})\leq P(\bigcap_{k}^{n}z_{k}\notin \pi(C_{k}))$. Thus, by \lemref{isoperimetric}, the claim follows.\\

Thirdly, let $(z(k))_{k}$ be an $\varepsilon-$ball walk that starts from a uniform point in $B_{L+n\varepsilon}(0)$. Also, let $z_{*}(k):=z(0)-z(k)$ so $z_{*}(k)$ is an $\varepsilon-$ball walk that starts from 0. Then uniform distribution of $z(0)$ yields
\begin{align*}
P(\forall k, z(k)\notin U_{k})&=P(z(0)\notin \bigcup_{k=0}^{n}(z_{*}(k)+U_{k}) )\\
&=1-\frac{E(vol(\bigcup_{k=0}^{n}(z_{*}(k)+U_{k})))}{vol(B_{L+n\varepsilon}(0))}.
\end{align*}
Therefore, 
\begin{align*}
1-\frac{E(vol(\bigcup_{k=0}^{n}(z_{*}(k)+\pi(C_{k})   )))}{vol(B_{L+n\varepsilon}(0))}&=P(\bigcap_{k}^{n}\pi(\zeta_{k})\notin \pi(C_{k}))\\
&\geq P(\forall k, z(k)\notin U_{k})-2P(\text{coupling fails})\\
&=1-\frac{E(vol(\bigcup_{k=0}^{n}(z_{*}(k)+U_{k})))}{vol(B_{L+n\varepsilon}(0))}-2P(\text{coupling fails})
\end{align*}

As $R\to \infty$ ,by \lemref{sphere}\ref{coupling} the $P(\text{coupling fails})=1-P(\forall k$ s.t. $\pi(\zeta(k))\neq z(k) )\to 0$. Also, the projection of the geodesic cap is a ball i.e. $\pi(C_{k})=B_{r_{k,R}}(0)\subset \mathbb{R}^{d}$. Let $r_{k}$ be radius of ball $U^{*}_{k}$. Then, by \lemref{sphere}\ref{Intermediate value property}
\begin{equation*} \lim \limits_{R\to \infty}r_{R,k}=r_{k} \end{equation*}

In other words, 
\begin{align*}
1-\frac{E(vol(\bigcup_{k=0}^{n}(z_{*}(k)+\pi(C_{k})   )))}{vol(B_{L+n\varepsilon}(0))}&=\lim \limits_{R\to \infty}1-\frac{E(vol(\bigcup_{k=0}^{n}(z_{*}(k)+B_{r_{k,R}}(0)   )))}{vol(B_{L+n\varepsilon}(0))}\\
&=\lim \limits_{R\to \infty}1-\frac{E(vol(\bigcup_{k=0}^{n}(z_{*}(k)+\pi(C_{k})   )))}{vol(B_{L+n\varepsilon}(0))}\\
&\geq \lim \limits_{R\to \infty}1-\frac{E(vol(\bigcup_{k=0}^{n}(z_{*}(k)+U_{k})))}{vol(B_{L+n\varepsilon}(0))}-2P(\text{coupling fails})\\
&=1-\frac{E(vol(\bigcup_{k=0}^{n}(z_{*}(k)+U_{k})))}{vol(B_{L+n\varepsilon}(0))}
\end{align*}

\end{proof}

\subsection{Research Problems}
\begin{enumerate}
\item  Does it also hold that for all $\lambda\in \mathbb{R}$ 
\begin{equation*}P[vol(\bigcup_{s\leq t}B_{s}+A_{s})\geq \lambda] \geq P[vol(\bigcup_{s\leq t}B_{s}+A^{*}_{s})] = E[vol(\bigcup_{s\leq t}S(B_{s},r_{s})\geq \lambda] .  \end{equation*}
\item Consider function $f:\mathbb{R}\to \mathbb{R}^{d}$. For which open sets A is the following quantity minized for $f\equiv 0$
\end{enumerate}
     \pagebreak

\pagebreak
\section{Spitzer's formula}
The electrostatic capacity of an object is defined by the following problem. Assume the object is conducting and charged so that its surface has a constant (unit) potential, and the potential outside the object decays to zero at infinite distance. The capacity can then be defined in terms of the asymptotic decay at large distances of the solution to the Dirichlet problem in the space surrounding the object. For dimension 2, this is is called the logarithmic capacity and for $n\geq 3$ the Newtonian capacity (for details see \cite{Landkoff_1972}). More specifically, the Newtonian capacity is defined as 
\begin{equation*}
Cap(A)=[inf_{\mu}\{\int_{A}\int_{A}G(x,y)d\mu(x)d\mu(y):\mu(A)=1~\text{probability measures}\}]^{-1},
\end{equation*}
where $G(x,y)=\frac{\Gamma(\frac{d}{2}-1)}{2\pi^{\frac{d}{2}}}\frac{1}{|x-y|^{d-2}}$ is the Green function. Spitzer proved that Cap(A) can be viewed as as the total heat $A\in \mathbb{R}^{3}$ can absorb \cite{Spitzer_1964}:
\begin{equation*}
Cap(A)=lim_{t\to\infty}\frac{1}{t}\int_{\mathbb{R}^{d}}P_{x}(\tau_{A}<t)dx=:lim_{t\to\infty}\frac{1}{t}E_{A}(t),
\end{equation*}
where $E_{A}(t)$ is called the energy of A. It was refined and extened to higher dimension by Le Gall \cite{le_gall_1990,Le_Gall_1988} and Port \cite{Port_1990}. In this section we will prove this formula.
     
\begin{proof}
First we prove that $\lim_{t\to \infty}E_{A}(t)-E_{A}(t-h)=h\cdot c$ for all $h>0$, where c will be later proved to be Cap(A) and thus for $h=t$ we obtain the result. Let $q_{A}(t,x,y)$ denote the probability of BM transitioning from $B_{0}=x$ to $B_{t}=y$ while conditioned to be killed when exiting A up to time t. In Port \cite{Port_Stone_1978} it is shown that $P_{x}(\tau_{A}>t, B_{t}\in B)=\int_{B}q_{A}(t,x,y)dy$ and $q_{A}(t,x,y)$ is symmetric in x,y. Using the Strong Markov property (SMP) and the above integral formula  in terms of $q_{A}(t,x,y)$, we obtain
\begin{align*}
E_{A}(t)-E_{A}(t-h)=&\int_{\mathbb{R}^{d}\setminus A}P_{x}(\tau_{A}\leq t)-P_{x}(\tau_{A}\leq t-h)dx\\
=&\int_{\mathbb{R}^{d}\setminus A}P_{x}(\tau_{A}\geq t-h)-P_{x}(\tau_{A}\geq t,B_{t}\in \mathbb{R}^{d} )dx\\
\stackrel{SMP}{=}&\int_{\mathbb{R}^{d}\setminus A}E_{x}[\int_{\mathbb{R}^{d}}p(\tau_{\mathbb{R}^{d}}-t+h,B_{\tau_{\mathbb{R}^{d}}},y)dy]dx-\int_{\mathbb{R}^{d}\setminus A}\int_{\mathbb{R}^{d}}q_{A}(t,x,y)dxdy.
\end{align*}
By definition of expectation 
\begin{equation*}
E_{x}[\int_{\mathbb{R}^{d}}p(\tau_{\mathbb{R}^{d}}-t+h,B_{\tau_{\mathbb{R}^{d}}},y)dy]=\int_{\mathbb{R}^{d}}\int_{\mathbb{R}^{d}}p(h,du,y)P_{x}(\tau_{\mathbb{R}^{d}}\geq t-h, B_{t}\in du)dy.
\end{equation*}
Thus, the symmetry of $q_{A}(t,x,y)$ gives
\begin{align*}
&\int_{\mathbb{R}^{d}\setminus A} \int_{\mathbb{R}^{d}}p(h,u,y)  \int_{\mathbb{R}^{d}} q_{A}(t-h,x,u) du dy     dx-\int_{\mathbb{R}^{d}\setminus A}q_{A}(t,x,y)dx\\
=& \int_{\mathbb{R}^{d}\setminus A} \int_{\mathbb{R}^{d}}p(h,y,u)  \int_{\mathbb{R}^{d}} q_{A}(t-h,u,x) du dy     dx-\int_{\mathbb{R}^{d}\setminus A}q_{A}(t,y,x)dx\\
=&\int_{\mathbb{R}^{d}}\int_{\mathbb{R}^{d}}p(h,y,u)P_{u}(\tau_{A}\geq t-h)du-P_{y}(\tau_{A}\geq t)dy\\
=&\int_{\mathbb{R}^{d}}P_{y}(h\leq \tau_{A^{c}}\leq t)-P_{y}(\tau_{A}\geq t)dy\\
=&\int_{\mathbb{R}^{d}}P_{y}(h\leq \tau_{A^{c}}\leq t, \tau_{A}\leq t)dy\\
=&\int_{\mathbb{R}^{d}}P_{y}(h\leq \tau_{A^{c}}\leq t, \tau_{A}\leq h)dy.
\end{align*}
Since the integrand $P_{y}(h\leq \tau_{A^{c}}\leq t, \tau_{A}\leq h)$ is monotone in t, monotone convergence yields
\begin{equation*}
lim_{t\to \infty}E_{A}(t)-E_{A}(t-h)=\int_{\mathbb{R}^{d}}P_{y}(\tau_{A}\leq h)dy.
\end{equation*}
Next since $\int_{\mathbb{R}^{d}}P_{y}(\tau_{A}\leq h)dy$ is an additive and monotone function of h, we get that $\int_{\mathbb{R}^{d}}P_{y}( \tau_{A}\leq h)dy=h\cdot c$. We will prove the stronger result $c=Cap(A)$, then for $h=t$ we get $lim_{t\to \infty}\frac{1}{t}E_{A}(t)=Cap(A)$.\\

Let $e_{h}(dy)=\psi_{h}(y)dy:=\frac{1}{h}P_{y}( \tau_{A}\leq h)dy$. We will show for the measures $e_{h}(dy)$ that
\begin{equation*}
c=lim_{h\to 0}\int e_{h}(dy):=lim_{h\to 0}\int\frac{1}{h}P_{y}( \tau_{A}\leq h)dy=\int e(dy)=Cap(A),
\end{equation*}
where $e(dy)$ is an equilibrium measure supported in A and corresponds to the potential $P_{y}(\tau_{A}<\infty)$ i.e. $P_{y}(\tau_{A}<\infty)=\int G(x,y)e(dy)$. Here are the steps:\\
Step 1: $P_{x}(\tau_{A}<\infty)=lim_{h\to 0}\int G(x,y)e_{h}(dy)=\int G(x,y)e(dy)$ and $e(\cdot)$ is supported on A.\\
Step 2: The minimizer for Cap(A) is achieved at $\frac{e(\cdot)}{e(A)}$ and this implies $lim_{h\to 0}\int\frac{1}{h}P_{y}( \tau_{A}\leq h)dy=\int e(dy)=e(A)=Cap(A)$.\\
\centerline{Step 1}
First we show that $P_{x}(\tau_{A}<\infty)=lim_{h\to 0}\int G(x,y)e_{h}(dy)$. Let $U_{h}$ be a uniform random variable on $[0,h]$, independent of BM and the stopping time $\tau_{A}$. Then for any bounded and continuous $f:\mathbb{R}^{d}\to \mathbb{R}$,
\begin{align*}
\int_{\mathbb{R}^{d}}f(y)G(x,y)\psi_{h}(y)dy=&\int_{0}^{\infty}p(t,x,y)f(y)\psi_{h}(y)dydt\\
=&\int_{0}^{\infty} E_{x}[f(B_{t})\cdot \psi_{h}(B_{t})]dt\\
=&\int_{0}^{\infty}\frac{1}{h}E_{x}[f(B_{t})1_{t\leq \tau_{A}\leq t+h}]dt\\
=&E_{x}[f(B_{\tau_{A}-U_{h}})1_{U_{h}<\tau_{A}}].
\end{align*}
Thus, $P_{x}(B_{\tau_{A}-U_{h}}\in B, U_{h}<\tau_{A})=\int_{B}G(x,y)\psi_{h}(y)dy$ for $B\subset A$ because indicator functions can be approximated by bounded continuous ones. By continuity of BM, $\lim_{h\to 0}\int_{B}G(x,y)\psi_{h}(y)dy=\lim_{h\to 0}P_{x}(B_{\tau_{A}-U_{h}}\in B, U_{h}<\tau_{A})=P_{x}[B_{\tau_{A}}\in B, \tau_{A}<\infty]$. Next we show that $lim_{h\to 0}\int G(x,y)e_{h}(dy)=\int G(x,y)e(dy)$. Because of continuity and boundedness of $y\mapsto \frac{1}{G(x,y)}$ on A we have 
\begin{align*}
lim_{h\to 0}e_{h}(dy)=&lim_{h\to 0}\psi_{h}(y)dy=lim_{h\to 0}\frac{1}{G(x,y)}P_{y}(B_{\tau_{A}-U_{h}}\in dy,U_{h}<\tau_{A})\\
=&\frac{1}{G(x,y)}P_{y}[B_{\tau_{A}}\in dy, \tau_{A}<\infty]dy=:e(dy).
\end{align*}
Thus, $P_{x}[B_{\tau_{A}}\in B, \tau_{A}<\infty]=\int_{B} G(x,y)e(dy)$ and in turn $e(A^{c})=0$ i.e. $e(\cdot)$ is supported on A.

\centerline{Step 2}
We will show that $\frac{e(\cdot)}{e(A)}$ is a minimizer for $\inf\limits_{\mu}\{\int\int G(x,y)d\mu(x)d\mu(y):\mu(A)=1~\text{probability measures}\}$. Let $I_{G}(\mu):=\int_{A}\int_{A}G(x,y)d\mu(x)d\mu(y)$. So for arbitrary measure $\mu$ on A with $\mu(A)=e(A)$, we will show that $I_{G}(\mu)\geq I_{G}(e)$. We first prove the following lemma

\begin{lem}
Let $\mu,\nu$ be finite measures on $\mathbb{R}^{d}$ and $\sigma:=\mu-\nu$, then
\begin{equation*}
\int_{\mathbb{R}^{d}}\int_{\mathbb{R}^{d}}G(x,y)d\sigma(x)d\sigma(y)\geq 0
~\text{and in the equality case, we obtain equality of measures }~\mu=\nu.
\end{equation*}
\end{lem}
\begin{proof}
From the semigroup property
\begin{equation*}
p(t,x,y)=\int p(\frac{t}{2}x,z)p(\frac{t}{2}z,y)dz.
\end{equation*}
Integrating wrt to $d\sigma(x)d\sigma(y)$ and using symmetry of $p(t,x,y)=p(t,y,x)$ yields
\begin{align*}
\int\int G(x,y)d\sigma(x)d\sigma(y)=&\int_{0}^{\infty}\int \int p(t,x,y)d\sigma(x)d\sigma(y)dt\\
=&\int_{0}^{\infty}\int \int \int p(\frac{t}{2}x,z)p(\frac{t}{2}z,y)dz d\sigma(x)d\sigma(y)dt\\
=&\int_{0}^{\infty}\int [\int p(\frac{t}{2}x,z)d\sigma(x)]^{2}dzdt\geq 0.
\end{align*}
Next we show equality of measures. Equality in the above inequality yields for a.e. z and t
\begin{equation*}
\int p(\frac{t}{2},x,z)d\sigma(x)=0.
\end{equation*}

Now fix a continuous function $f:\mathbb{R}^{d}\to [0,\infty)$ with compact support. It suffices to prove $\int f(x)d\sigma(x)=0$ because indicator functions can be approximated by bounded continuous ones. By the delta-function property of transition probability $p(t,x,y)$, we obtain
\begin{equation*}
f(x)=lim_{t\to 0}\int f(z)p(\frac{t}{2},x,z)dz.
\end{equation*}
Thus, 
\begin{equation*}
\int f(x)d\sigma(x)=\int lim_{t\to 0}\int f(z)p(\frac{t}{2},x,z)dz d\sigma(x)=lim_{t\to 0}\int f(z) \int p(\frac{t}{2},x,z)dz d\sigma(x)dz =0.
\end{equation*}
\end{proof}
Now we show that $I_{G}(\mu)\geq I_{G}(e)$ where $\mu(A)=e(A)$. Since $P_{x}(\tau_{A}<\infty)=1$ for $x\in A$, it yields
\begin{equation*}
\int\int G(x,y)e(dx)e(dy)=\int_{A}P_{x}(\tau_{A}<\infty)e(dy)=e(A).
\end{equation*}
Therefore, by above lemma 
\begin{align*}
I_{G}(\mu)-I_{G}(e)=&I_{G}(\mu)-e(A)\\
\geq &I_{G}(\mu)+e(A)-2\int_{A}P_{x}(\tau_{A}<\infty)d\mu(y)\\
=&I_{G}(\mu)+I_{G}(e)-2\int \int G(x,y)e(dx)d\mu(y)\\
=&\int \int G(x,y) d(e-\mu)(x)d(e-\mu)(y)\geq 0\\
\Rightarrow & I_{G}(\mu)\geq I_{G}(e).
\end{align*}
Thus, probability measure $\frac{e(\cdot)}{e(A)}$ is a minimizer. In fact the lemma shows that this is the unique minimizer: if $I_{G}(\mu)=I_{G}(e)\Rightarrow \int \int G(x,y) d(e-\mu)(x)d(e-\mu)(y)=0\Rightarrow e(\cdot)=\mu(\cdot)$. Therefore, we showed that 
\begin{align*}
Cap(A)=&[inf_{\mu}\{I_{G}(\mu):\mu(A)=1~\text{probability measures}\}]^{-1}\\
=&[\int\int G(x,y)\frac{e(dx)}{e(A)}\frac{e(dy)}{e(A)}]^{-1}\\
=&e(A)[\int_{A} P_{x}(\tau_{A}<\infty)\frac{e(dy)}{e(A)}]^{-1}\\
=&e(A)=lim_{h\to 0}\int\frac{1}{h}P_{y}(\tau_{A}\leq h)dy\\
=&lim_{t\to \infty}\frac{1}{h}(E_{A}(t)-E_{A}(t-h))\\
=&lim_{t\to \infty}\frac{1}{t}E_{A}(t) ~ \text{for h=t}.
\end{align*}

\end{proof}

\subsection{Research problems}
Le Gall, Port \cite{le_gall_1990,Le_Gall_1988,Port_1990} improved the estimate of $E_{A}(t)$ for compact $A\in \mathbb{R}^{d}$ to \cite{Le_Gall_1988}:
\begin{align*}
\text{For}~ d=3\\
E_{A}(t)=&Cap(A)t+\frac{4}{(2\pi)^{\frac{3}{2}}}Cap(A)^{2}t^{\frac{1}{2}}+\frac{1}{2\pi^{2}}Cap(A)^{3}-|A|-\frac{1}{2\pi}\int \int |z-y|e(dz)e(dy)]+O(t^{\frac{-1}{2}}).\\
\text{For}~ d=4\\
E_{A}(t)=&Cap(A)t+\frac{1}{(2\pi)^{2}}Cap(A)^{2}logt+\frac{1}{(2\pi)^{2}}[(1+log2-\gamma)Cap(A)^{2}-2\int \int log|z-y|e(dz)e(dy)]\\
&+\frac{1}{8\pi^{4}}Cap(A)^{3}\frac{logt}{t}-|A|+o(\frac{logt}{t}).\\
\text{For}~ d\geq 5\\
E_{A}(t)=&Cap(A)t+\frac{\Gamma(\frac{d}{2}-2)}{4\pi^{\frac{d}{2}}}\int \int|z-y|^{4-d}e(dz)e(dy)-|A|-\frac{4}{(2\pi)^{\frac{d}{2}}(d-2)(d-4)}Cap(A)^{2}t^{2-\frac{d}{2}}\\
&+O(t^{1-\frac{d}{2}}).
\end{align*}
We denoted the Euler constant by $\gamma$. For latest results see Van de Berg \cite{Van_Berg_2007}. 
                
\textbf{Acknowledgements} \\               
It is a pleasure to thank Almut Burchard for introducing me to the topic and for valuable discussions. 
    
\newcommand{\etalchar}[1]{$^{#1}$}

\end{document}